# GENERAL MAXIMUM LIKELIHOOD EMPIRICAL BAYES ESTIMATION OF NORMAL MEANS

BY WENHUA JIANG AND CUN-HUI ZHANG[1]

*Rutgers University*

We propose a general maximum likelihood empirical Bayes (GMLEB) method for the estimation of a mean vector based on observations with i.i.d. normal errors. We prove that under mild moment conditions on the unknown means, the average mean squared error (MSE) of the GMLEB is within an infinitesimal fraction of the minimum average MSE among all separable estimators which use a single deterministic estimating function on individual observations, provided that the risk is of greater order than $(\log n)^5/n$. We also prove that the GMLEB is uniformly approximately minimax in regular and weak $\ell_p$ balls when the order of the length-normalized norm of the unknown means is between $(\log n)^{\kappa_1}/n^{1/(p\wedge 2)}$ and $n/(\log n)^{\kappa_2}$. Simulation experiments demonstrate that the GMLEB outperforms the James–Stein and several state-of-the-art threshold estimators in a wide range of settings without much down side.

**1. Introduction.** This paper concerns the estimation of a vector with i.i.d. normal errors under the average squared loss. The problem, known as the compound estimation of normal means, has been considered as the canonical model or motivating example in the developments of empirical Bayes, admissibility, adaptive nonparametric regression, variable selection, multiple testing and many other areas in statistics. It also carries significant practical relevance in statistical applications since the observed data are often understood, represented or summarized as the sum of a signal vector and the white noise.

There are three main approaches in the compound estimation of normal means. The first one is general empirical Bayes (EB) [27, 30], which assumes

Received December 2007; revised July 2008.
[1]Supported by NSF Grants DMS-05-04387 and DMS-06-04571 and NSA Grant MDS 904-02-1-0063.
*AMS 2000 subject classifications.* 62C12, 62G05, 62G08, 62G20, 62C25.
*Key words and phrases.* Compound estimation, empirical Bayes, adaptive estimation, white noise, shrinkage estimator, threshold estimator.

essentially no knowledge about the unknown means but still aims to attain the performance of the oracle separable estimator based on the knowledge of the empirical distribution of the unknowns. Here a separable estimator is one that uses a fixed deterministic function of the $i$th observation to estimate the $i$th mean. This greedy approach, also called nonparametric EB [26], was proposed the earliest among the three, but it is also the least understood, in spite of [28, 29, 30, 36, 37, 38]. Efron [15] attributed this situation to the lack of applications with many unknowns before the information era and pointed out that "current scientific trends favor a greatly increased role for empirical Bayes methods" due to the prevalence of large, high-dimensional data and the rapid rise of computing power. The methodological and theoretical challenge, which we focus on in this paper, is to find the "best" general EB estimators and sort out the type and size of problems suitable for them.

The second approach, conceived with the celebrated Stein's proof of the inadmissibility of the optimal unbiased estimator and the introduction of the James–Stein estimator [22, 31], is best understood through its parametric or linear EB interpretations [16, 17, 26]. The James–Stein estimator is minimax over the entire space of the unknown mean vector and well approximates the optimal linear separable estimator based on the oracular knowledge of the first two empirical moments of the unknown means. Thus, it achieves the general EB optimality when the empirical distribution of the unknown means are approximately normal. However, the James–Stein estimator does not perform well by design compared with the general EB when the minimum risk of linear separable estimators is far different from that of all separable estimators [36]. Still, what is the cost of being greedy with the general EB when the empirical distribution of the unknown means is indeed approximately normal?

The third approach focuses on unknown mean vectors which are sparse in the sense of having many (near) zeros. Such sparse vectors can be treated as members of small $\ell_p$ balls with $p < 2$. Examples include the estimation of functions with unknown discontinuity or inhomogeneous smoothness across different parts of a domain in nonparametric regression or density problems [13]. For sparse means, the James–Stein or the oracle linear estimators could perform much worse than threshold estimators [12]. Many threshold methods have been proposed and proved to possess (near) optimality properties for sparse signals, including the universal [13], SURE [14], FDR [1, 2], the generalized $C_p$ [3] and the parametric EB posterior median (EBThresh) [24]. These estimators can be viewed as approximations of the optimal candidate in certain families of separable threshold estimators, so that they do not perform well by design compared with the general EB when the minimum risk of separable threshold estimators is far different from that of all separable estimators [38]. Again, what is the cost of being greedy with the general EB when the unknown means are indeed very sparse?

Since general EB methods have to spend more "degrees of freedom" for nonparametric estimation of its oracle rule, compared with linear and threshold methods, the heart of the question is whether the gain by aiming at the smaller general EB benchmark risk is large enough to offset the additional cost of the nonparametric estimation.

We propose a general maximum likelihood EB (GMLEB) in which we first estimate the empirical distribution of the unknown means by the generalized maximum likelihood estimator (MLE) [25] and then plug the estimator into the oracle general EB rule. In other words, we treat the unknown means as i.i.d. variables with a completely unknown common "prior" distribution (for the purpose of deriving the GMLEB, whether the unknowns are actually deterministic or random), estimate the nominal prior with the generalized MLE, and then use the Bayes rule for the estimated prior. The basic idea was discussed in the last paragraph of [27] as a general way of deriving solutions to compound decision problems, although the notion of MLE was vague at that time without a parametric model and not much has been done since then about using the generalized MLE to estimate the nominal prior in compound estimation.

Our results affirm that by aiming at the minimum risk of *all* separable estimators, the greedier general EB approach realizes significant risk reduction over linear and threshold methods for a wide range of the unknown signal vectors for moderate and large samples, and this is especially so for the GMLEB. We prove that the risk of the GMLEB estimator is within an infinitesimal fraction of the general EB benchmark when the risk is of the order $n^{-1}(\log n)^5$ or greater depending on the magnitude of the weak $\ell_p$ norm of the unknown means, $0 < p \leq \infty$. Such adaptive ratio optimality is obtained through a general oracle inequality which also implies the adaptive minimaxity of the GMLEB over a broad collection of regular and weak $\ell_p$ balls. This adaptive minimaxity result unifies and improves upon the adaptive minimaxity of threshold estimators for sparse means [1, 14, 24] and the Fourier general EB estimators for moderately sparse and dense means [38]. We demonstrate the superb risk performance of the GMLEB for moderate samples through simulation experiments, and describe algorithms to show its computational feasibility.

The paper is organized as follows. In Section 2, we highlight our results and formally introduce certain necessary terminologies and concepts. In Section 3 we provide upper bounds for the regret of a regularized Bayes rule using a predetermined and possibly misspecified prior. In Section 4 we prove an oracle inequality for the GMLEB, compared with the general EB benchmark risk. The consequences of this oracle inequality, including statements of our adaptive ratio optimality and adaptive minimaxity results in full strength, are also discussed in Section 4. In Section 5 we present more simulation results. Section 6 contains some discussion. Mathematical proofs of theorems, propositions and lemmas are given either right after their statements or in the Appendix.

**2. Problem formulation and highlight of main results.** Let $X_i$ be independent statistics with

$$X_i \sim \varphi(x-\theta_i) \sim N(\theta_i, 1), \qquad i=1,\dots,n, \tag{2.1}$$

under a probability measure $P_{n,\boldsymbol{\theta}}$, where $\boldsymbol{\theta}=(\theta_1,\dots,\theta_n)$ is an unknown signal vector. Our problem is to estimate $\boldsymbol{\theta}$ under the compound loss

$$L_n(\widehat{\boldsymbol{\theta}}, \boldsymbol{\theta}) = n^{-1}\|\widehat{\boldsymbol{\theta}}-\boldsymbol{\theta}\|^2 = \frac{1}{n}\sum_{i=1}^{n}(\widehat{\theta}_i-\theta_i)^2 \tag{2.2}$$

for any given estimator $\widehat{\boldsymbol{\theta}}=(\widehat{\theta}_1,\dots,\widehat{\theta}_n)$. Throughout this paper, the unknown means $\theta_i$ are assumed to be deterministic as in the standard compound decision theory [27]. To avoid confusion, the Greek $\theta$ is used only with boldface as a deterministic mean vector $\boldsymbol{\theta}$ in $\mathbb{R}^n$ or with subscripts as elements of $\boldsymbol{\theta}$. A random mean is denoted by $\xi$ as in (2.3) below. The estimation of i.i.d. random means is discussed in Section 6.3.

We divide the section into seven subsections to describe (1) the general and restricted EB, (2) the GMLEB method, (3) the computation of the GMLEB, (4) some simulation results, (5) the adaptive ratio optimality of the GMLEB, (6) the adaptive minimaxity of the GMLEB in $\ell_p$ balls and (7) minimax theory in $\ell_p$ balls.

Throughout the paper, boldface letters denote vectors and matrices, for example, $\mathbf{X}=(X_1,\dots,X_n)$, $\varphi(x)=e^{-x^2/2}/\sqrt{2\pi}$ denotes the standard normal density, $\widetilde{L}(y)=\sqrt{-\log(2\pi y^2)}$ denotes the inverse of $y=\varphi(x)$ for positive $x$ and $y$, $x\vee y=\max(x,y)$, $x\wedge y=\min(x,y)$, $x_+=x\vee 0$ and $a_n\asymp b_n$ means $0<a_n/b_n+b_n/a_n=O(1)$. In a number of instances, $\log(x)$ should be viewed as $\log(x\vee e)$. Univariate functions are applied to vectors per component. Thus, an estimator of $\boldsymbol{\theta}$ is separable if it is of the form $\widehat{\boldsymbol{\theta}}=t(\mathbf{X})=(t(X_1),\dots,t(X_n))$ with a predetermined Borel function $t(\cdot)$. In the vector notation, it is convenient to state (2.1) as $\mathbf{X}\sim N(\boldsymbol{\theta},\mathbf{I}_n)$ with $\mathbf{I}_n$ being the identity matrix in $\mathbb{R}^n$.

2.1. *Empirical Bayes.* The compound estimation of a vector of determinist normal means is closely related to the Bayes estimation of a single random mean. In this Bayes problem, we estimate a univariate random parameter $\xi$ based on a univariate $Y$ such that

$$Y|\xi \sim N(\xi,1), \qquad \xi\sim G, \qquad \text{under } P_G. \tag{2.3}$$

The prior distribution $G = G_n$ which naturally matches the unknown means $\{\theta_i, i \le n\}$ in (2.1) is the empirical distribution

$$G_n(u) = G_{n,\boldsymbol{\theta}}(u) = \frac{1}{n}\sum_{i=1}^{n} I\{\theta_i \le u\}. \tag{2.4}$$

Here and in the sequel, subscripts $_{n,\boldsymbol{\theta}}$ indicate dependence of distribution or probability upon $n$ and the unknown deterministic vector $\boldsymbol{\theta}$.

The *fundamental theorem of compound decisions* [27] in the context of the $\ell_2$ loss asserts that the compound risk of a separable rule $\widehat{\boldsymbol{\theta}} = t(\mathbf{X})$ under the probability $P_{n,\boldsymbol{\theta}}$ in the multivariate model (2.1) is identical to the MSE of the same rule $\widehat{\xi} = t(Y)$ under the prior (2.4) in the univariate model (2.3):

$$E_{n,\boldsymbol{\theta}} L_n(t(\mathbf{X}), \boldsymbol{\theta}) = E_{G_n}(t(Y) - \xi)^2. \tag{2.5}$$

For any true or nominal priors $G$, denote the Bayes rule as

$$t_G^* = \arg\min_t E_G(t(Y) - \xi)^2 = \frac{\int u\varphi(Y-u)G(du)}{\int \varphi(Y-u)G(du)} \tag{2.6}$$

and the minimum Bayes risk as

$$R^*(G) = E_G(t_G^*(Y) - \xi)^2, \tag{2.7}$$

where the minimum is taken over all Borel functions. It follows from (2.5) that among all separable rules, the compound risk is minimized by the Bayes rule with prior (2.4), resulting in the general EB benchmark

$$R^*(G_n) = E_{n,\boldsymbol{\theta}} L_n(t_{G_n}^*(\mathbf{X}), \boldsymbol{\theta}) = \min_{t(\cdot)} E_{n,\boldsymbol{\theta}} L_n(t(\mathbf{X}), \boldsymbol{\theta}). \tag{2.8}$$

The general EB approach seeks procedures which approximate the Bayes rule $t_{G_n}^*(\mathbf{X})$ or approximately achieve the risk benchmark $R^*(G_n)$ in (2.8).

Given a class of functions $\mathscr{D}$, the aim of the restricted EB is to attain

$$R_{\mathscr{D}}(G_n) = \inf_{t\in\mathscr{D}} E_{n,\boldsymbol{\theta}} L_n(t(\mathbf{X}), \boldsymbol{\theta}) = \inf_{t\in\mathscr{D}} E_{G_n}(t(Y) - \xi)^2, \tag{2.9}$$

approximately. This provides EB interpretations for all the adaptive methods discussed in the Introduction, with $\mathscr{D}$ being the classes of all linear functions for the James–Stein estimator, all soft threshold functions for the SURE [14], and all hard threshold functions for the generalized $C_p$ [3] or the FDR [1]. For the EBThresh [24], $\mathscr{D}$ is the class of all posterior median functions $t(y) = \text{median}(\xi|Y = y)$ under the probability $P_G$ in (2.3) for priors of the form

$$G(u) = \omega_0 I\{0 \le u\} + (1 - \omega_0)G_0(u/\tau), \tag{2.10}$$

where $\omega_0$ and $\tau$ are free and $G_0$ is given [e.g., $dG_0(u)/du = e^{-|u|}/2$].

Compared with linear and threshold methods, the general EB approach is greedier since it aims at the smaller benchmark risk: $R^*(G_n) \le R_{\mathscr{D}}(G_n)$ for all $\mathscr{D}$. This could still backfire when the regret

$$r_{n,\boldsymbol{\theta}}(\widehat{t}_n) = E_{n,\boldsymbol{\theta}} L_n(\widehat{t}_n(\mathbf{X}), \boldsymbol{\theta}) - R^*(G_n) \tag{2.11}$$

of using an estimator $\widehat{t}_n(\cdot)$ of the general EB oracle rule $t^*_{G_n}(\cdot)$ is greater than the difference $R_{\mathscr{D}}(G_n) - R^*(G_n)$ in benchmarks, but our simulation and oracle inequalities prove that $r_{n,\boldsymbol{\theta}}(\widehat{t}_n) = o(1)R^*(G_n)$ uniformly for a wide range of the unknown vector $\boldsymbol{\theta}$ and moderate/large samples.

Zhang [36] proposed a general EB method based on a Fourier infinite-order smoothing kernel. The Fourier general EB estimator is asymptotically minimax over the entire parameter space and approximately reaches the general EB benchmark (2.8) uniformly for dense and moderately sparse signals, provided that the oracle Bayes risk is of the order $n^{-1/2}(\log n)^{3/2}$ or greater [36]. Hybrid general EB estimators have been developed [38] to combine the features and optimality properties of the Fourier general EB and threshold estimators. Still, the performance of general EB methods is sometimes perceived as uncertain in moderate samples [24]. Indeed, the Fourier general EB requires selection of certain tuning parameters and its proven theoretical properties are not completely satisfying. This motivates our investigation.

2.2. *The GMLEB.* The GMLEB method replaces the unknown prior $G_n$ of the oracle rule $t^*_{G_n}$ by its generalized MLE [25]

$$\widehat{G}_n = \widehat{G}_n(\cdot; \mathbf{X}) = \mathop{\arg\max}_{G \in \mathscr{G}} \prod_{i=1}^n f_G(X_i), \tag{2.12}$$

where $\mathscr{G}$ is the family of *all* distribution functions and $f_G$ is the density

$$f_G(x) = \int \varphi(x-u) G(du) \tag{2.13}$$

of the normal location mixture by distribution $G$.

The estimator (2.12) is called the generalized MLE since the likelihood is used only as a vehicle to generate the estimator. The $G$ here is used only as a nominal prior. In our adaptive ratio and minimax optimality theorems and oracle inequality, the GMLEB is evaluated under the measures $P_{n,\boldsymbol{\theta}}$ in (2.1) where the unknowns $\theta_i$ are assumed to be deterministic parameters.

Since (2.12) is typically solved by iterative algorithms, we allow approximate solutions to be used. For definiteness and notation simplicity, the generalized MLE in the sequel is any solution of

$$\widehat{G}_n \in \mathscr{G}, \qquad \prod_{i=1}^n f_{\widehat{G}_n}(X_i) \ge q_n \sup_{G \in \mathscr{G}} \prod_{i=1}^n f_G(X_i) \tag{2.14}$$

with $q_n = (e\sqrt{2\pi}/n^2) \wedge 1$, although the theoretical results in this paper all hold verbatim for less stringent (2.14) with $0 \le \log(1/q_n) \le c_0(\log n)$ for any fixed constant $c_0$. Formally, the GMLEB estimator is defined as

$$\widehat{\boldsymbol{\theta}} = t^*_{\widehat{G}_n}(\mathbf{X}) \quad \text{or equivalently} \quad \widehat{\theta}_i = t^*_{\widehat{G}_n}(X_i), \qquad i = 1, \ldots, n, \tag{2.15}$$

where $t^*_G$ is the Bayes rule in (2.6) and $\widehat{G}_n$ is any approximate generalized MLE (2.14) for the nominal prior (2.4). Clearly, the GMLEB estimator (2.15) is completely nonparametric and does not require any restriction, regularization, bandwidth selection or other forms of tuning.

The GMLEB is location equivariant in the sense that

$$t^*_{\widehat{G}_n(\cdot;\mathbf{X}+c\mathbf{e})}(\mathbf{X} + c\mathbf{e}) = t^*_{\widehat{G}_n(\cdot;\mathbf{X})}(\mathbf{X}) + c\mathbf{e} \tag{2.16}$$

for all real $c$, where $\mathbf{e} = (1, \ldots, 1) \in \mathbb{R}^n$. This is due to the location equivariance of the generalized MLE: $\widehat{G}_n(x; \mathbf{X} + c\mathbf{e}) = \widehat{G}_n(x - c; \mathbf{X})$. Compared with the Fourier general EB estimators [36, 38], the GMLEB (2.15) is more appealing since the function $t^*_{\widehat{G}_n}(x)$ of $x$ enjoys all analytical properties of Bayes rules: monotonicity, infinite differentiability and more. However, the GMLEB is much harder to analyze than the Fourier general EB. We first address the computational issues in the next section.

2.3. *Computation of the GMLEB.* It follows from the Carathéodory's theorem [9] that there exists a discrete solution of (2.12) with no more than $n + 1$ support points. A discrete approximate generalized MLE $\widehat{G}_n$ with $m$ support points can be written as

$$\widehat{G}_n = \sum_{j=1}^m \widehat{w}_j \delta_{u_j}, \qquad \widehat{w}_j \ge 0, \qquad \sum_{j=1}^m \widehat{w}_j = 1, \tag{2.17}$$

where $\delta_u$ is the probability distribution giving its entire mass to $u$. Given (2.17), the GMLEB estimator can be easily computed as

$$\widehat{\theta}_i = t^*_{\widehat{G}_n}(X_i) = \frac{\sum_{j=1}^m u_j \varphi(X_i - u_j)\widehat{w}_j}{\sum_{j=1}^m \varphi(X_i - u_j)\widehat{w}_j}, \tag{2.18}$$

since $t^*_G(x)$ is the conditional expectation as in (2.6).

Since the generalized MLE $\widehat{G}_n$ is completely nonparametric, the support points $\{u_j, j \le m\}$ and weights $\{\widehat{w}_j, j \le m\}$ in (2.17) are selected or computed solely to maximize the likelihood in (2.12). There are quite a few possible algorithms for solving (2.14), but all depend on iterative approximations. Due to the monotonicity of $\varphi(t)$ in $t^2$, the generalized MLE (2.12) puts all its mass in the interval $I_0 = [\min_{1\le i\le n} X_i, \max_{1\le i\le n} X_i]$. Given a fine grid $\{u_j\}$ in $I_0$, the EM-algorithm [11, 35]

$$\widehat{w}_j^{(k)} = \frac{1}{n}\sum_{i=1}^n \frac{\widehat{w}_j^{(k-1)}\varphi(X_i - u_j)}{\sum_{\ell=1}^m \widehat{w}_\ell^{(k-1)}\varphi(X_i - u_\ell)} \tag{2.19}$$

optimizes the weights $\{\widehat{w}_j\}$. In Section 6.2, we provide a conservative statistical criterion on $\{u_j\}$ and an EM-stopping rule to guarantee (2.14).

We took a simple approach in our simulation experiments. Given $\{X_i, 1 \le i \le n\}$ and with $X_0 = 0$, we chose the grid points $\{u_j\}$ as a set of multipliers of $\varepsilon = \max_{0\le i<j\le n} |X_i - X_j|/999$ with $u_j = u_{j-1} + \varepsilon$ and the range

$$-j_0\varepsilon = u_1 - \varepsilon < \min_{0\le i\le n} X_i \le u_1, \qquad u_m = (m - j_0)\varepsilon \le \max_{0\le i\le n} X_i < u_m + \varepsilon$$

with an integer $j_0 \in [1, m]$. This ensures $u_{j_0} = 0$ as a grid point and $999 \le m \le 1000$. We ran 100 EM-iterations (2.19) in our simulations. We have tried to optimize both the support points $\{u_j\}$ and weights $\{\widehat{w}_j\}$ in the EM-algorithm, but gained limited improvements.

The GMLEB estimator (2.18) depends slightly on the initialization of the EM-algorithm due to the nonuniqueness of the GMLEB estimator and the fixed number of EM-iterations in our implementation. Since the generalized MLE (2.12) is unique only up to the values of $\{f_{\widehat{G}_n}(X_i), i \le n\}$, different EM-initializations lead to different versions of $\widehat{G}_n$, which then result in different values of $t^*_{\widehat{G}_n}(X_i)$ in (2.18). This nonuniqueness persists even when we run infinitely many EM-iterations. Nevertheless, our theoretical results hold for all versions of the GMLEB.

We consider two options in our simulation experiments. The first option initializes the weights with the uniform distribution $\widehat{w}_j = 1/m$. The second option takes into consideration of the possible sparsity of the signal by putting a good starting mass at $u_{j_0} = 0$:

$$\widehat{w}_{j_0} = \widehat{\omega}_0, \qquad \widehat{w}_j = \frac{1 - \widehat{\omega}_0}{m - 1}, \qquad j \ne j_0. \tag{2.20}$$

We estimate the proportion of zeros within the $n$ means by a Fourier method,

$$\widehat{\omega}_0 = \frac{1}{n}\sum_{j=1}^{n} \psi(X_j; h_n), \qquad \psi(z; h) = \int h\psi_0(ht)e^{t^2/2}\cos(zt)\,dt$$

as in [32, 33], where $\psi_0$ is a density function with support $[-1, 1]$ and $h_n = \{\kappa(\log n)\}^{-1/2}$ is the bandwidth, $\kappa \le 1$. In our simulation experiments, the uniform $[-1, 1]$ density is used as $\psi_0$ and $\kappa = 1/2$. To distinguish the two options of initializing the EM-algorithm, we reserve the name GMLEB for the uniform initialization and call (sparse-) S-GMLEB the estimator with the initialization (2.20) when we report simulation results.

2.4. *Some simulation results.* Johnstone and Silverman [24] reported results of an extensive simulation study of 18 threshold estimators, including eight options of their EBThresh, the SURE and adaptive SURE [14], the FDR [1] at three levels, three block threshold methods [7, 8] and the soft and hard threshold at the universal threshold level $\sqrt{2\log n}$. In their simulations, the overall best performer is the EBThresh using the posterior median for the prior (2.10) with the double exponential $dG_0(u)/du = e^{-|u|}/2$ and the MLE of $(\omega_0, \tau)$.

TABLE 1
*Average total squared errors $\|\widehat{\boldsymbol{\theta}}-\boldsymbol{\theta}\|^2$ for $n=1000$ unknown means in various binary models where $\theta_j$ is either 0 or $\mu$ with the number of nonzero $\theta_i=\mu$ being 5, 50 or 500. The "Best" stands for the best simulation results in Table 1 of Johnstone and Silverman [24].*
*Each entry is based on 100 replications*

| **# nonzero** | **5** | | | | **50** | | | | **500** | | | |
|---|---|---|---|---|---|---|---|---|---|---|---|---|
| $\boldsymbol{\mu}$ | **3** | **4** | **5** | **7** | **3** | **4** | **5** | **7** | **3** | **4** | **5** | **7** |
| James–Stein | 45 | 76 | 113 | 199 | 312 | 442 | 556 | 716 | 822 | 889 | 933 | 954 |
| EBThresh | **37** | **34** | **20** | **8** | **212** | **151** | **103** | 74 | 862 | 873 | 792 | 653 |
| SURE | 42 | 64 | 73 | 75 | 416 | 609 | 215 | 214 | 835 | 834 | 842 | 828 |
| FDR (0.01) | 43 | 54 | 29 | **6** | 388 | 299 | 132 | **57** | 2587 | 1322 | 667 | 520 |
| FDR (0.1) | 42 | 38 | **21** | 13 | 278 | 163 | 115 | 99 | 1162 | 744 | 662 | 640 |
| GMLEB | **39** | **34** | 23 | 11 | **157** | **105** | **58** | **14** | **459** | **285** | **139** | **18** |
| S-GMLEB | **32** | **28** | **17** | **6** | **150** | **99** | **54** | **10** | **454** | **282** | **136** | **15** |
| F-GEB | 94 | 94 | 89 | 88 | 223 | 185 | 135 | 103 | **520** | **363** | **237** | **131** |
| HF-GEB | 37 | 34 | 20 | 8 | 197 | 150 | 99 | 72 | 499 | 334 | 192 | 83 |
| "Best" | 34 | 32 | 17 | 5 | 201 | 156 | 95 | 52 | 829 | 730 | 609 | 505 |
| Oracle | 27 | 22 | 12 | 0.8 | 144 | 93 | 46 | 3 | 443 | 273 | 128 | 8 |

In Table 1, we display our simulation results under exactly the same setting as in [24] for nine estimators: the James–Stein, the EBThresh [24] using the double exponential $dG_0$ in (2.10) and the MLE of $(\omega_0, \tau)$, the SURE [14], the FDR [1] at levels $q = 0.01$ and $q = 0.1$, the GMLEB (2.15) with the uniform initialization, the S-GMLEB with the initialization (2.20), the F-GEB and HF-GEB as the Fourier general EB [36] and a hybrid [38] of its monotone version with the EBThresh. In each column, boldface entries denote the top three performers other than the hybrid estimator. We also display as "Best" the best of the simulation results in [24] over the 18 threshold estimators and as Oracle the average simulated risk of the oracle Bayes rule $t^*_{G_n}$ in (2.8).

These simulation results can be summarized as follows. The average $\ell_2$ loss of the S-GMLEB happens to be the smallest among the nine estimators, with the S-GMLEB and GMLEB clearly outperforming all other methods by large margins for dense and moderately sparse signals. For very sparse signals, the S-GMLEB, the EBThresh, the GMLEB and the FDR estimators yield comparable results, and they all outperform the Fourier general EB and James–Stein estimators. Compared with the oracle, the regrets of the S-GMLEB and GMLEB are nearly fixed constants. Since the oracle prior (2.4) has a point mass at 0 in all the models used to generate data in this simulation experiment, the S-GMLEB yields slightly better results than the GMLEB as expected. The hybrid estimator correctly switches to the EBThresh for very sparse signals.

These simulations and more presented in Section 5 demonstrate the computational affordability of the proposed GMLEB. The most surprising aspect of the results in Table 1 is the strong performance of the both versions of the GMLEB for the most sparse signals with 0.5% of $\theta_i$ being nonzero, since the GMLEB is not specially designed to recover such signals (and threshold estimators are).

2.5. *Adaptive ratio optimality.* Our theoretical results match well with the supreme performance of the GMLEB in our simulation experiments. We describe here the adaptive ratio optimality of the GMLEB and in the next section the adaptive minimaxity of the GMLEB in $\ell_p$ balls.

The adaptive ratio optimality holds for an estimator $\widehat{\boldsymbol{\theta}}:\mathbf{X}\to\mathbb{R}^n$ if its risk is uniformly within a fraction of the general EB benchmark

$$\sup_{\boldsymbol{\theta}\in\Theta_n^*}\frac{E_{n,\boldsymbol{\theta}}L_n(\widehat{\boldsymbol{\theta}},\boldsymbol{\theta})}{R^*(G_{n,\boldsymbol{\theta}})}\le 1+o(1) \tag{2.21}$$

in certain classes $\Theta_n^*\subset\mathbb{R}^n$ of the unknown vector $\boldsymbol{\theta}$, where $L_n(\cdot,\cdot)$ is the average squared loss (2.2), $G_{n,\boldsymbol{\theta}}=G_n$ is the empirical distribution of the unknowns in (2.4) and $R^*(G_n)$ is the general EB benchmark risk (2.8) achieved by the oracle Bayes rule $t^*_{G_n}(\mathbf{X})$.

THEOREM 1. *Let $\mathbf{X}\sim N(\boldsymbol{\theta},\mathbf{I}_n)$ under $P_{n,\boldsymbol{\theta}}$ with a deterministic $\boldsymbol{\theta}\in\mathbb{R}^n$. Let $t^*_{\widehat{G}_n}(\cdot)$ be the GMLEB in (2.15) with an approximate solution $\widehat{G}_n$ satisfying (2.14). Let $G_n=G_{n,\boldsymbol{\theta}}$ and $R^*(G)$ be as in (2.4) and (2.7). Then,*

$$\frac{E_{n,\boldsymbol{\theta}}L_n(t^*_{\widehat{G}_n}(\mathbf{X}),\boldsymbol{\theta})}{R^*(G_n)}=\frac{E_{n,\boldsymbol{\theta}}\|t^*_{\widehat{G}_n}(\mathbf{X})-\boldsymbol{\theta}\|^2}{\min_t E_{n,\boldsymbol{\theta}}\|t(\mathbf{X})-\boldsymbol{\theta}\|^2}\le 1+o(1) \tag{2.22}$$

*for the compound loss (2.2), provided that for certain constants $b_n$*

$$\frac{nR^*(G_n)}{(\sqrt{\log n}\vee\max_{i\le n}|\theta_i-b_n|)(\log n)^{9/2}}\to\infty.$$

*In particular, if $\max_{i\le n}|\theta_i-b_n|=O(\sqrt{\log n})$ and $nR^*(G_n)/(\log n)^5\to\infty$, then (2.22) holds.*

For any sequences of constants $M_n\to\infty$, Theorem 1 provides the adaptive ratio optimality (2.21) of the GMLEB in the classes

$$\Theta_n^*=\{\boldsymbol{\theta}\in\mathbb{R}^n:R^*(G_{n,\boldsymbol{\theta}})\ge M_n n^{-1}(\log n)^{9/2}(\sqrt{\log n}\vee\|\boldsymbol{\theta}\|_\infty)\}.$$

This is a consequence of an oracle inequality for the GMLEB $\widehat{t}_n = t^*_{\widehat{G}_n}$ in Section 4.2, which uniformly bound from the above

$$\widetilde{r}_{n,\boldsymbol{\theta}}(\widehat{t}_n) = \sqrt{E_{n,\boldsymbol{\theta}} L_n(\widehat{t}_n(\mathbf{X}), \boldsymbol{\theta})} - \sqrt{R^*(G_n)} \tag{2.23}$$

in terms of the weak $\ell_p$ norm of $\boldsymbol{\theta}$. The quantity (2.23) can be viewed as the regret for the minimization of the square root of the MSE, instead of (2.11). Clearly, $r_{n,\boldsymbol{\theta}}(\widehat{t}_n)/R^*(G_n) \le o(1)$ iff $\widetilde{r}_{n,\boldsymbol{\theta}}(t^*_{\widehat{G}_n})/\sqrt{R^*(G_n)} \le o(1)$. A more general version of Theorem 1 is given in Section 4.3.

In the EB literature, the asymptotic optimality of $\widehat{\boldsymbol{\theta}}$ is defined as

$$G_n \xrightarrow{\mathrm{D}} G \quad \Rightarrow \quad E_{n,\boldsymbol{\theta}} L_n(\widehat{\boldsymbol{\theta}}, \boldsymbol{\theta}) - R^*(G_n) \to 0 \tag{2.24}$$

for deterministic vectors $\boldsymbol{\theta} \in \mathbb{R}^n$ [27, 36]. In the EB model

$$(Y_i, \xi_i) \text{ i.i.d.}, \qquad Y_i | \xi_i \sim N(\xi_i, 1), \qquad \xi_i \sim G, \qquad \text{under } P_G \tag{2.25}$$

with data $\{Y_i\}$, the EB asymptotic optimality is defined as

$$\lim_{n\to\infty} E_G \sum_{i=1}^n (\widehat{\xi}_i - \xi_i)^2 / n = R^*(G). \tag{2.26}$$

We call (2.21) adaptive ratio optimality since it is much stronger than both notions of asymptotic optimality in its uniformity in $\boldsymbol{\theta} \in \Theta^*_n$ and its focus on the harder standard of the relative error, due to $R^*(G_n) \le E_{n,\boldsymbol{\theta}} L_n(\mathbf{X}, \boldsymbol{\theta}) = 1$. The difference among these optimality properties is significant for moderate samples in view of some very small $R^*(G_n) \approx \text{Oracle}/1000$ in Table 1.

Theorem 1 is location invariant, since the GMLEB is location equivariant by (2.16) and $R^*(G_n)$ is location invariant by (2.8). Thus, if $\theta_i = b_n$ for most $i \le n$, the GMLEB performs equally well whether $b_n = 0$ or not. Moreover, if $\theta_i \in B\ \forall i$ for a finite set $B \subset \mathbb{R}$, the GMLEB adaptively shrinks toward the points in $B$ [19]. This is evident in Table 1 for $\#\{i : \theta_i = 7\} \in \{50, 500\}$ with $B = \{0, 7\}$. In fact, if $\#\{x : x \in B_n\} = O(1)$ and $\min_{B_n \ni x \ne y \in B_n} |x - y| \to \infty$, then $G_n(B_n) = 1$ implies $R^*(G_n) \to 0$. Threshold methods certainly do not possess these location invariance and multiple shrinkage properties.

2.6. *Adaptive minimaxity in $\ell_p$ balls.* Minimaxity is commonly used to measure the performance of statistical procedures. For $\Theta \subset \mathbb{R}^n$, the minimax risk for the average squared loss (2.2) is

$$\mathscr{R}_n(\Theta) = \inf_{\widetilde{\boldsymbol{\theta}}} \sup_{\boldsymbol{\theta} \in \Theta} E_{n,\boldsymbol{\theta}} L_n(\widetilde{\boldsymbol{\theta}}, \boldsymbol{\theta}), \tag{2.27}$$

where the infimum is taken over all Borel mappings $\widetilde{\boldsymbol{\theta}} : \mathbf{X} \to \mathbb{R}^n$. An estimator is minimax in a specific class $\Theta$ of unknown mean vectors if it attains $\mathscr{R}_n(\Theta)$, but this does not guarantee satisfactory performance since the minimax estimator is typically uniquely tuned to the specific set $\Theta$. For small $\Theta$, the minimax estimator has high risk outside $\Theta$. For large $\Theta$, the minimax estimator is too conservative by focusing on the worst case scenario within $\Theta$. Adaptive minimaxity overcomes this difficulty by requiring

$$\frac{\sup_{\boldsymbol{\theta}\in\Theta_n} E_{n,\boldsymbol{\theta}} L_n(\widehat{\boldsymbol{\theta}},\boldsymbol{\theta})}{\mathscr{R}_n(\Theta_n)} \to 1 \tag{2.28}$$

uniformly for a wide range of sequences $\{\Theta_n \subset \mathbb{R}^n, n \ge 1\}$ of parameter classes. Define (regular or strong) $\ell_p$ balls as

$$\Theta_{p,C,n} = \left\{\boldsymbol{\theta} = (\theta_1,\dots,\theta_n) : n^{-1}\sum_{i=1}^n |\theta_i|^p \le C^p\right\}. \tag{2.29}$$

The quantity $C$ in (2.29), called length-normalized or standardized radius of the $\ell_p$ ball, is denoted as $\eta$ in [1, 12, 24], where adaptive minimaxity in $\ell_p$ balls with $C = C_n \to 0$ and $p < 2$ is used to measure the performance of estimators for sparse $\boldsymbol{\theta}$. The following theorem establishes the adaptive minimaxity of the GMLEB in $\ell_p$ balls with radii $C = C_n$ in intervals diverging to $(0,\infty)$. This covers sparse and dense $\boldsymbol{\theta}$ simultaneously. Adaptive minimaxity of the GMLEB in weak $\ell_p$ balls is discussed in Section 4.3.

THEOREM 2. *Let $\mathbf{X} \sim N(\boldsymbol{\theta}, \mathbf{I}_n)$ under $P_{n,\boldsymbol{\theta}}$ with a deterministic $\boldsymbol{\theta} \in \mathbb{R}^n$. Let $\widehat{\boldsymbol{\theta}} = t^*_{\widehat{G}_n}(\mathbf{X})$ be the GMLEB in (2.15) with an approximate solution $\widehat{G}_n$ satisfying (2.14). Let $L_n(\cdot,\cdot)$ be the average squared loss (2.2) and $\mathscr{R}_n(\Theta)$ be the minimax risk (2.27). Then, as $n \to \infty$, the adaptive minimaxity (2.28) holds in $\ell_p$ balls (2.29) with $\Theta_n = \Theta_{p,C_n,n}$, provided that*

$$\frac{n^{1/(p\wedge 2)} C_n}{(\log n)^{\kappa_1(p)}} \to \infty, \qquad \frac{C_n}{n}(\log n)^{\kappa_2(p)} \to 0, \tag{2.30}$$

*where $\kappa_1(p) = 1/2 + 4/p + 3/p^2$ for $p < 2$, $\kappa_1(2) = 13/4$, $\kappa_1(p) = 5/2$ for $p > 2$, and $\kappa_2(p) = 9/2 + 4/p$.*

Theorem 2 is a consequence of the oracle inequality in Section 4.2 and the minimax theory in [12]. An outline of this argument is given in the next section. An alternative statement of the conclusion of Theorem 2 is

$$\lim_{(n,M)\to(\infty,\infty)} \sup_{C\in\mathscr{C}_{p,n}(M)} \frac{\sup_{\boldsymbol{\theta}\in\Theta_{p,C,n}} E_{n,\boldsymbol{\theta}} L_n(t^*_{\widehat{G}_n}(\mathbf{X}),\boldsymbol{\theta})}{\mathscr{R}_n(\Theta_{p,C,n})} = 1,$$

where $\mathscr{C}_{p,n}(M) = [Mn^{-1/(p\wedge 2)}(\log n)^{\kappa_1(p)}, n/\{M(\log n)^{\kappa_2(p)}\}]$. In Section 4.3, we offer an analogues result for weak $\ell_p$ balls. The powers $\kappa_1(p)$ and $\kappa_2(p)$ of the logarithmic factors in (2.30) and in the definition of $\mathscr{C}_{p,n}(M)$ are crude. This is further discussed in Section 6.

Adaptive and approximate minimax estimators of the normal means in $\ell_p$ balls have been considered in [1, 3, 12, 14, 24, 36, 38]. Donoho and Johnstone [14] proved that as $(n, C_n) \to (\infty, 0+)$, with $nC_n^p/(\log n)^{p/2} \to \infty$ for $p < 2$,

$$\mathscr{R}_n(\Theta_{p,C_n,n}) = (1 + o(1)) \min_{t \in \mathscr{D}} \max_{\boldsymbol{\theta} \in \Theta_{p,C_n,n}} E_{n,\boldsymbol{\theta}} L_n(t(\mathbf{X}), \boldsymbol{\theta}), \tag{2.31}$$

where $\mathscr{D}$ is the collection of all (soft or hard) threshold rules. Therefore, adaptive minimaxity (2.28) in small $\ell_p$ balls $\Theta_n = \Theta_{p,C_n,n}$ can be achieved by threshold rules with suitable data-driven threshold levels. This has been done using the FDR [1] for $(\log n)^5/n \le C_n^p \le n^{-\kappa}$ with $p < 2$ and any $\kappa > 0$. Zhang [38] proved that (2.28) holds for the Fourier general EB estimator of [36] in $\Theta_n = \Theta_{p,C_n,n}$ for $C_n^p\sqrt{n}/(\log n)^{1+(p\wedge 2)/2} \to \infty$.

A number of estimators have been proven to possess the adaptive rate minimaxity in the sense of attaining within a bounded factor of the minimax risk. In $\ell_p$ balls $\Theta_{p,C_n,n}$, the EBThresh is adaptive rate minimax for $p \le 2$ and $nC_n^p \ge (\log n)^2$ [24], while the generalized $C_p$ is adaptive rate minimax for $p < 2$ and $1 \le O(1)nC_n^p$ [3]. It follows from [3, 38] that a hybrid between the Fourier general EB and universal soft threshold estimators is also adaptive rate minimax in $\Theta_{p,C_n,n}$ for $1 \le O(1)nC_n^p$.

The adaptive minimaxity as provided in Theorem 2 unifies the adaptive minimaxity of different types estimators in different ranges of the radii $C_n$ of the $\ell_p$ balls with the exception of the two very extreme ends, due to the crude power $\kappa_1(p)$ of the logarithmic factor for small $C_n$ and the requirement of an upper bound for large $C_n$. The hybrid Fourier general EB estimator achieves the adaptive rate minimaxity in a wider range of $\ell_p$ balls than what we prove here for the GMLEB. However, as we have seen in Table 1, the finite sample performance of the GMLEB is much stronger. It seems that the less stringent and commonly considered adaptive rate minimaxity leaves too much room to provide adequate indication of finite sample performance.

2.7. *Minimax theory in $\ell_p$ balls.* Instead of the general EB approach, adaptive minimax estimation in small $\ell_p$ balls can be achieved by threshold methods, provided that the radius is not too small. However, since (2.31) does not hold for fixed $p > 0$ and $C \in (0, \infty)$, threshold estimators are not asymptotically minimax with $\Theta_n = \Theta_{p,C,n}$ in (2.28) for fixed $(p, C)$. Consequently, adaptive minimax estimations in small, fixed and large $\ell_p$ balls are often treated separately in the literature. In this section, we explain the general EB approach for adaptive minimax estimation, which provides a unified treatment for $\ell_p$ balls of different ranges of radii. This provides an outline of the proof of Theorem 2. Minimax theory in weak $\ell_p$ balls will be discussed in Section 4.3.

We first discuss the relationship between the minimax estimation of a deterministic vector $\boldsymbol{\theta}$ in $\ell_p$ balls and the minimax estimation of a single random mean under an unknown "prior" in $L_p$ balls. For positive $p$ and $C$, the $L_p$ balls of distribution functions are defined as

$$\mathscr{G}_{p,C} = \left\{ G : \int |u|^p G(du) \le C^p \right\}.$$

Since $\mathscr{G}_{p,C}$ is a convex class of distributions, the minimax theorem provides

$$\mathscr{R}(\mathscr{G}_{p,C}) = \min_t \max_{G \in \mathscr{G}_{p,C}} E_G(t(Y) - \xi)^2 = \max_{G \in \mathscr{G}_{p,C}} R^*(G) \le 1 \tag{2.32}$$

for the estimation of a single real random parameter $\xi$ in the model (2.3), where $R^*(G)$ is the minimum Bayes risk in (2.7). Thus, since $G_n = G_{n,\boldsymbol{\theta}} \in \mathscr{G}_{p,C}$ for $\boldsymbol{\theta} \in \Theta_{p,C,n}$, the fundamental theorem of compound decisions (2.5) implies that (2.32) dominates the compound minimax risk (2.27) in $\ell_p$ balls:

$$\mathscr{R}_n(\Theta_{p,C,n}) \le \inf_{t(x)} \sup_{\boldsymbol{\theta} \in \Theta_{p,C,n}} E_{n,\boldsymbol{\theta}} L_n(t(\mathbf{X}), \boldsymbol{\theta}) \le \mathscr{R}(\mathscr{G}_{p,C}) \le 1. \tag{2.33}$$

Donoho and Johnstone [12] proved that as $C^{p\wedge 2} \to 0+$

$$\left| \frac{\mathscr{R}(\mathscr{G}_{p,C})}{C^{p\wedge 2}\{2\log(1/C^p)\}^{(1-p/2)_+}} - 1 \right| \to 0 \tag{2.34}$$

and that for either $p \ge 2$ with $C_n > 0$ or $p < 2$ with $nC_n^p/(\log n)^{p/2} \to \infty$,

$$\left| \frac{\mathscr{R}_n(\Theta_{p,C_n,n})}{\mathscr{R}(\mathscr{G}_{p,C_n})} - 1 \right| \to 0. \tag{2.35}$$

In the general EB approach, the aim is to find an estimator $\widehat{t}_n$ of $t^*_{G_n}$ with small regret (2.11) or (2.23). If the approximation to $t^*_{G_n}$ in risk is sufficiently accurate and uniformly within a small fraction of $\mathscr{R}(\mathscr{G}_{p,C_n})$ for $\boldsymbol{\theta} \in \Theta_{p,C_n,n}$, the maximum risk of the general EB estimator in $\Theta_{p,C_n,n}$ would be within the same small fraction of $\mathscr{R}(\mathscr{G}_{p,C_n})$, since the risk of $t^*_{G_n}$ is bounded by $R^*(G_{n,\boldsymbol{\theta}}) \le \mathscr{R}(\mathscr{G}_{p,C_n})$ for $\boldsymbol{\theta} \in \Theta_{p,C_n,n}$. Thus, (2.35) plays a crucial role in general EB.

It follows from (2.23), (2.32) and (2.29) that

$$\sup_{\boldsymbol{\theta} \in \Theta_{p,C,n}} \sqrt{E_{n,\boldsymbol{\theta}} L_n(\widehat{t}_n(\mathbf{X}), \boldsymbol{\theta})} \le \sup_{\boldsymbol{\theta} \in \Theta_{p,C,n}} \widetilde{r}_{n,\boldsymbol{\theta}}(\widehat{t}_n) + \sqrt{\mathscr{R}(\mathscr{G}_{p,C})}. \tag{2.36}$$

Thus, by (2.34) and (2.35), the adaptive minimaxity (2.28) of $\widehat{\boldsymbol{\theta}} = \widehat{t}_n(\mathbf{X})$ in $\ell_p$ balls $\Theta_n = \Theta_{p,C_n,n}$ is a consequence of an oracle inequality of the form

$$\sup_{\boldsymbol{\theta} \in \Theta_{p,C_n,n}} \widetilde{r}_{n,\boldsymbol{\theta}}(\widehat{t}_n) = o(1)\sqrt{J_{p,C_n}} \tag{2.37}$$

with $J_{p,C} = \min\{1, C^{p\wedge 2}\{1 \vee (2\log(1/C^p))\}^{(1-p/2)_+}\}$. In our proof, (2.34) and the upper bound $\mathscr{R}(\mathscr{G}_{p,C}) \le 1$ provide $\inf_C \mathscr{R}(\mathscr{G}_{p,C})/J_{p,C} > 0$. Although $J_{p,C}$ provides the order of $\mathscr{R}(\mathscr{G}_{p,C})$ for each $p$ via (2.34), explicit expressions of the minimax risk $\mathscr{R}_n(\Theta_{p,C,n})$ for general fixed $(p, C, n)$ or the minimax risk $\mathscr{R}(\mathscr{G}_{p,C})$ for fixed $(p, C)$ with $p \ne 2$ are still open problems.

**3. A regularized Bayes estimator with a misspecified prior.** In this section, we consider a fixed probability $P_{G_0}$ under which

$$Y|\xi \sim N(\xi, 1), \qquad \xi \sim G_0. \tag{3.1}$$

Recall [5, 28] that for the estimation of a normal mean, the Bayes rule (2.6) and its risk (2.7) can be expressed in terms of the mixture density $f_G(x)$ as

$$t_G^*(x) = x + \frac{f_G'(x)}{f_G(x)}, \qquad R^*(G) = 1 - \int \left(\frac{f_G'}{f_G}\right)^2 f_G, \tag{3.2}$$

in the model (2.3), where $f_G(x) = \int \varphi(x-u)G(du)$ is as in (2.13).

Suppose the true prior $G_0$ is unknown but a deterministic approximation of it, say $G$, is available. The Bayes formula (3.2) could still be used, but we may want to avoid dividing by a near-zero quantity. This leads to the following regularized Bayes estimator:

$$t_G^*(x;\rho) = x + \frac{f_G'(x)}{f_G(x) \vee \rho}. \tag{3.3}$$

For $\rho = 0$, $t_G^*(x;0) = t_G^*(x)$ is the Bayes estimator for the prior $G$. For $\rho = \infty$, $t_G^*(x;\infty) = x$ gives the MLE of $\xi$ which requires no knowledge of the prior. The following proposition, proved in the Appendix, describes some analytical properties of the regularized Bayes estimator.

PROPOSITION 1. *Let $\widetilde{L}(y) = \sqrt{-\log(2\pi y^2)}$, $y \ge 0$, be the inverse function of $y = \varphi(x)$. Then, the value of the regularized Bayes estimator $t_G^*(x;\rho)$ in (3.3) is always between those of the Bayes estimator $t_G^*(x)$ in (2.6) and the MLE $t_G^*(x;\infty) = x$. Moreover, for all real $x$*

$$\begin{cases} |x - t_G^*(x;\rho)| = \dfrac{|f_G'(x)|}{f_G(x) \vee \rho} \le \widetilde{L}(\rho), & 0 < \rho \le (2\pi e)^{-1/2}, \\ 0 \le (\partial/\partial x) t_G^*(x;\rho) \le \widetilde{L}^2(\rho), & 0 < \rho \le (2\pi e^3)^{-1/2}. \end{cases} \tag{3.4}$$

REMARK 1. In [36], a slightly different inequality

$$\left(\frac{f_G'(x)}{f_G(x)}\right)^2 \frac{f_G(x)}{f_G(x) \vee \rho} \le \widetilde{L}^2(\rho), \qquad 0 \le \rho < (2\pi e^2)^{-1/2}, \tag{3.5}$$

was used to derive oracle inequalities for Fourier general EB estimators. The extension to the derivative of $t_G^*(x;\rho)$ here is needed for the application of the Gaussian isoperimetric inequality in Proposition 4.

The next theorem provides oracle inequalities which bound the regret of using (3.3) due to the lack of the knowledge of the true $G_0$. Let

$$d(f,g) = \left(\int (f^{1/2} - g^{1/2})^2\right)^{1/2} \tag{3.6}$$

denote the Hellinger distance. The upper bounds assert that the regret is no greater than the square of the Hellinger distance between the mixture densities $f_G$ and $f_{G_0}$ up to certain logarithmic factors.

THEOREM 3. *Suppose (3.1) holds under $P_{G_0}$. Let $t^*_G(x;\rho)$ be the regularized Bayes rule in (3.3) with $0<\rho\le(2\pi e^2)^{-1/2}$. Let $f_G$ be as in (2.13).*

(i) *There exists a universal constant $M_0$ such that*

$$\begin{aligned}&[E_{G_0}\{t^*_G(Y;\rho)-\xi\}^2-R^*(G_0)]^{1/2}\\ &\quad\le M_0\max\{|\log\rho|^{3/2},|\log(d(f_G,f_{G_0}))|^{1/2}\}d(f_G,f_{G_0})\\ &\qquad+\left\{\int\left(1-\frac{f_{G_0}}{\rho}\right)^2_+\frac{(f'_{G_0})^2}{f_{G_0}}\right\}^{1/2},\end{aligned}\tag{3.7}$$

*where $R^*(G_0)=E_{G_0}\{t^*_{G_0}(Y)-\xi\}^2$ is the minimum Bayes risk in (2.7).*

(ii) *If $\int_{|u|>x_0}G_0(du)\le M_1|\log\rho|^3\varepsilon_0^2$ and $2(x_0+1)\rho\le M_2|\log\rho|^2\varepsilon_0^2$ for a certain $\varepsilon_0\ge d(f_G,f_{G_0})$ and finite positive constants $\{x_0,M_1,M_2\}$, then*

$$\begin{aligned}&E_{G_0}\{t^*_G(Y;\rho)-\xi\}^2-R^*(G_0)\\ &\quad\le 2(M_0+M_1+M_2)\max(|\log\rho|^3,|\log\varepsilon_0|)\varepsilon_0^2,\end{aligned}\tag{3.8}$$

*where $M_0$ is a universal constant.*

REMARK 2. For $G=G_0$ (3.7) becomes an identity, so that the square of the first term on the right-hand side of (3.7) represents an upper bound for the regret of using a misspecified $G$ in the regularized Bayes estimator (3.3) instead of the true $G_0$ for the same regularization level $\rho$. Under the additional tail probability condition on $G_0$ and for sufficiently small $\rho$, (3.8) provides an upper bound for the regret of not knowing $G_0$, compared with the Bayes estimator (3.2) with the true $G=G_0$.

REMARK 3. Since the second term on the right-hand side of (3.7) is increasing in $\rho$ and the first is logarithmic in $1/\rho$, we are allowed to take $\rho>0$ of much smaller order than $d(f_G,f_{G_0})$ in (3.7), for example, under moment conditions on $G_0$. Still, the cubic power of the logarithmic factors in (3.7) and (3.8) is crude.

The following lemma plays a crucial role in the proof of Theorem 3.

LEMMA 1. *Let $d(f,g)$ be as in (3.6) and $\widetilde{L}(y)=\sqrt{-\log(2\pi y^2)}$. Then,*

$$\int\frac{(f'_G-f'_{G_0})^2}{f_G\vee\rho+f_{G_0}\vee\rho}\le e^2 2d^2(f_G,f_{G_0})\max(\widetilde{L}^6(\rho),2a^2)\tag{3.9}$$

*for $\rho\le 1/\sqrt{2\pi}$, where $a^2=\max\{\widetilde{L}^2(\rho)+1,|\log d^2(f_G,f_{G_0})|\}$.*

PROOF OF THEOREM 3. Let

$$\|g\|_h = \left\{ \int g^2(x)h(x)\,dx \right\}^{1/2}$$

be the $L_2(h(x)\,dx)$ norm for $h \geq 0$. Since $t^*_{G_0}$ is the Bayes rule, by (3.3)

$$\begin{aligned}
&[E_{G_0}\{t^*_G(Y;\rho) - \xi\}^2 - E_{G_0}\{t^*_{G_0}(Y) - \xi\}^2]^{1/2} \\
&\quad = \|f'_G/(f_G \vee \rho) - f'_{G_0}/f_{G_0}\|_{f_{G_0}} \qquad (3.10)\\
&\quad \leq r(f_G,\rho) + \|(1 - f_{G_0}/\rho)_+ f'_{G_0}/f_{G_0}\|_{f_{G_0}},
\end{aligned}$$

where $r(f_G,\rho) = \|f'_G/(f_G \vee \rho) - f'_{G_0}/(f_{G_0} \vee \rho)\|_{f_{G_0}}$.

Let $w_* = 1/(f_G \vee \rho + f_{G_0} \vee \rho)$. For $G_1 = G$ or $G_1 = G_0$,

$$\begin{aligned}
\int \left( \frac{f'_{G_1}}{f_{G_1} \vee \rho} - 2f'_{G_1} w_* \right)^2 f_{G_0} &\leq \int \left( \frac{f'_{G_1}}{f_{G_1} \vee \rho} |f_G - f_{G_0}| w_* \right)^2 f_{G_0} \\
&\leq \widetilde{L}^2(\rho) \int (f_G - f_{G_0})^2 w_*^2 f_{G_0}
\end{aligned}$$

due to $|f'_{G_1}|/(f_{G_1} \vee \rho) \leq \widetilde{L}(\rho)$ by (3.4). Since $(\sqrt{f_G} + \sqrt{f_{G_0}})^2 w_* \leq 2$ and $w_* f_{G_0} \leq 1$, we find

$$\begin{aligned}
r(f_G,\rho) &\leq 2\|(f'_G - f'_{G_0})w_*\|_{f_{G_0}} + 2\widetilde{L}(\rho)\|(f_G - f_{G_0})w_*\|_{f_{G_0}} \\
&\leq 2\|f'_G - f'_{G_0}\|_{w_*} + 2\widetilde{L}(\rho)\sqrt{2}d(f_G, f_{G_0}).
\end{aligned}$$

Thus, (3.7) follows from (3.10) and (3.9).

To prove (3.8) we use Lemma 6.1 in [38]:

$$\begin{aligned}
&\int_{f_{G_0} < \rho} \left( \frac{f'_{G_0}}{f_{G_0}} \right)^2 f_{G_0} \\
&\quad \leq \int_{|u| > x_0} G_0(du) + 2x_0\rho \max\{\widetilde{L}^2(\rho), 2\} + 2\rho\sqrt{\widetilde{L}^2(\rho) + 2} \\
&\quad \leq (M_1 + M_2)|\log \rho|^3 \varepsilon_0^2,
\end{aligned}$$

due to $|\log \rho| \geq \widetilde{L}^2(\rho) \geq 2$. This and (3.7) imply (3.8). $\square$

**4. An oracle inequality for the GMLEB.** In this section, we provide an oracle inequality which bound the regret (2.23) and thus (2.11) of using the GMLEB $t^*_{\widehat{G}_n}$ in (2.15) against the oracle Bayes rule $t^*_{G_n}$ in (2.8). We provide the main elements leading to the oracle inequality in Section 4.1 before presenting the oracle inequality and an outline of its proof in Section 4.2. Section 4.3 discusses the consequences of the oracle inequality, including a sharper version of Theorem 1 and the adaptive minimaxity in weak $\ell_p$ balls.

4.1. *Elements leading to the oracle inequality.* It follows from the fundamental theorem of compound decisions (2.5) that for separable estimators $\widehat{\boldsymbol{\theta}} = t(\mathbf{X})$, the compound risk is identical to the MSE of $\widehat{\xi} = t(Y)$ for the estimation of a single real random parameter $\xi$ under $P_G$ in (2.3), so that Theorem 3 provides an upper bound for the regret of the regularized Bayes rule $t^*_G(\mathbf{X};\rho)$ in terms of the Hellinger distance $d(f_G, f_{G_n})$ and $\rho > 0$. We have proved in [39] a large deviation upper bound for the Hellinger distance $d(f_{\widehat{G}_n}, f_{G_n})$. We will show that the GMLEB estimator $t^*_{\widehat{G}_n}(\mathbf{X})$ is identical to its regularized version $t^*_{\widehat{G}_n}(\mathbf{X};\rho_n)$ for certain $|\log\rho_n| \asymp \log n$ when the generalized MLE (2.12) or its approximation (2.14) are used. Still, $t^*_{\widehat{G}_n}(\mathbf{X};\rho_n)$ is not separable, since the generalized MLE $\widehat{G}_n$ is based on the same data $\mathbf{X}$. A natural approach of deriving oracle inequalities is then to combine Theorem 3 with a maximal inequality. This requires in addition an entropy bound for the class of regularized Bayes rules $t^*_G(x;\rho)$ with given $\rho > 0$ and an exponential inequality for the difference between the loss and risk for each regularized Bayes rule. In the rest of this section, we provide these crucial components of our theoretical investigation.

4.1.1. *A large deviation inequality for the convergence of an approximate generalized MLE.* Under the i.i.d. assumption of the EB model (2.25), Ghosal and van der Vaart [20] obtained an exponential inequality for the Hellinger loss of the generalized MLE of a normal mixture density in terms of the $L_\infty$ norm of $\theta_i$. This result can be improved upon using their newer entropy calculation in [21]. The results in [20, 21] are unified and further improved upon in the i.i.d. case and extended to deterministic $\boldsymbol{\theta} = (\theta_1, \dots, \theta_n)$ in weak $\ell_p$ balls for all $0 < p \le \infty$ in [39]. This latest result, stated below as Theorem 4, will be used here in conjunction of Theorem 3 to prove oracle inequalities for the GMLEB.

The $p$th weak moment of a distribution $G$ is

$$\mu^w_p(G) = \left\{\sup_{x>0} x^p \int_{|u|>x} G(du)\right\}^{1/p} \tag{4.1}$$

with $\mu^w_\infty(G) = \inf\{x : \int_{|u|>x} G(du) = 0\}$. Define convergence rates

$$\begin{aligned}\varepsilon(n,G,p) &= \max[\sqrt{2\log n}, \{n^{1/p}\sqrt{\log n}\,\mu^w_p(G)\}^{p/(2+2p)}]\sqrt{\frac{\log n}{n}} \\ &= \max\left[\sqrt{\frac{2\log n}{n}}, \left\{\sqrt{\log n}\frac{\mu^w_p(G)}{n}\right\}^{p/(2+2p)}\right]\sqrt{\log n}\end{aligned} \tag{4.2}$$

with $\varepsilon(n,G,\infty) = \{(2\log n) \vee (\sqrt{\log n}\,\mu^w_\infty(G))\}^{1/2}\sqrt{(\log n)/n}$.

THEOREM 4. *Let* $\mathbf{X} \sim N(\boldsymbol{\theta}, \mathbf{I}_n)$ *under* $P_{n,\boldsymbol{\theta}}$ *with a deterministic* $\boldsymbol{\theta} \in \mathbb{R}^n$*. Let* $f_G$ *and* $G_n$ *be as in (2.13) and (2.4), respectively. Let* $\widehat{G}_n$ *be a certain approximate generalized MLE satisfying (2.14). Then, there exists a universal constant* $x_*$ *such that for all* $x \ge x_*$ *and* $\log n \ge 2/p$*,*

$$P_{n,\boldsymbol{\theta}}\{d(f_{\widehat{G}_n}, f_{G_n}) \ge x\varepsilon_n\} \le \exp\left(-\frac{x^2 n\varepsilon_n^2}{2\log n}\right) \le e^{-x^2 \log n}, \tag{4.3}$$

*where* $\varepsilon_n = \varepsilon(n, G_n, p)$ *is as in (4.2) and* $d(f,g)$ *is the Hellinger distance (3.6). In particular, for any sequences of constants* $M_n \to \infty$ *and fixed positive* $\alpha$ *and* $c$*,*

$$\varepsilon_n \asymp \begin{cases} n^{-p/(2+2p)}(\log n)^{(2+3p)/(4+4p)}, \\ \qquad \text{if } \mu_p^w(G_n) = O(1) \text{ with a fixed } p, \\ n^{-1/2}(\log n)^{3/4}\{M_n^{1/2} \vee (\log n)^{1/4}\}, \\ \qquad \text{if } G_n([-M_n, M_n]) = 1 \text{ and } p = \infty, \\ n^{-1/2}(\log n)^{1/(2(2\wedge\alpha))+3/4}, \\ \qquad \text{if } \int e^{|cu|^\alpha} G_n(du) = O(1) \text{ and } p \asymp \log n. \end{cases}$$

REMARK 4. Under the condition $G([-M_n, M_n]) = 1$ and the i.i.d. assumption (2.25) with $G$ depending on $n$, the large deviation bound in [20] provides the convergence rate $\varepsilon_n \asymp n^{-1/2}(\log n)^{1/2}\{M_n \vee (\log n)^{1/2}\}$, and the entropy calculation in [21] leads to the convergence rate $\varepsilon_n \asymp n^{-1/2}(\log n)\sqrt{M_n}$. These rates are slower than the rate in Theorem 4 when $M_n/\sqrt{\log n} \to \infty$.

REMARK 5. The proof of Theorem 4 is identical for the generalized MLE (2.12) and its approximation (2.14). The constant $x_*$ is universal for $q_n = (e\sqrt{2\pi}/n^2) \wedge 1$ in (2.14) and depends on $\sup_n |\log q_n|/\log n$ in general.

4.1.2. *Representation of the GMLEB estimator as a regularized one at data points.* The connection between the GMLEB estimator (2.15) and the regularized Bayes rule (3.3) in Theorem 3 is provided by

$$t^*_{\widehat{G}_n}(\mathbf{X}) = t^*_{\widehat{G}_n}(\mathbf{X}; \rho_n), \qquad \rho_n = q_n/(en\sqrt{2\pi}), \tag{4.4}$$

where $q_n$ is as in (2.14). This is a consequence of the following proposition.

PROPOSITION 2. *Let* $f(x|u)$ *be a given family of densities and* $\{X_i, i \le n\}$ *be given data. Let* $\widehat{G}_n$ *be an approximate generalized MLE of a mixing distribution satisfying*

$$\prod_{i=1}^n \int f(X_i|u)\widehat{G}_n(du) \ge q_n \sup_G \prod_{i=1}^n \int f(X_i|u)G(du)$$

*for certain* $0 < q_n \le 1$. *Then, for all* $j = 1, \ldots, n$

$$f_{\widehat{G}_n}(X_j) = \int f(X_j|u)\widehat{G}_n(du) \ge \frac{q_n}{en} \sup_u f(X_j|u).$$

*In particular, (4.4) holds for* $f(x|u) = \varphi(x-u)$.

PROOF. Let $j$ be fixed and $u_j = \arg\max f(X_j|u)$. Define $\widehat{G}_{n,j} = (1-\varepsilon)\widehat{G}_n + \varepsilon\delta_{u_j}$ with $\varepsilon = 1/n$, where $\delta_u$ is the unit mass at $u$. Since $f(x|u) \ge 0$, $f_{\widehat{G}_{n,j}}(X_i) \ge (1-\varepsilon) f_{\widehat{G}_n}(X_i)$ and $f_{\widehat{G}_{n,j}}(X_j) \ge \varepsilon f(X_j|u_j)$, so that

$$\frac{1}{q_n}\prod_{i=1}^n f_{\widehat{G}_n}(X_i) \ge \prod_{i=1}^n f_{\widehat{G}_{n,j}}(X_i) \ge (1-\varepsilon)^{n-1}\varepsilon f(X_j|u_j)\prod_{i\ne j} f_{\widehat{G}_n}(X_i).$$

Thus, $f_{\widehat{G}_n}(X_j) \ge q_n(1-\varepsilon)^{n-1}\varepsilon f(X_j|u_j)$ with $\varepsilon = 1/n$, after the cancellation of $f_{\widehat{G}}(X_i)$ for $i \ne j$. The conclusion follows from $(1-1/n)^{n-1} \ge 1/e$. $\square$

4.1.3. *An entropy bound for regularized Bayes rules.* We now provide an entropy bound for collections of regularized Bayes rules. For any family $\mathscr{H}$ of functions and semidistance $d_0$, the $\varepsilon$-covering number is

$$N(\varepsilon, \mathscr{H}, d_0) = \inf\left\{N : \mathscr{H} \subseteq \bigcup_{j=1}^N \mathrm{Ball}(h_j, \varepsilon, d_0)\right\} \tag{4.5}$$

with $\mathrm{Ball}(h, \varepsilon, d_0) = \{f : d_0(f, h) < \varepsilon\}$. For each fixed $\rho > 0$ define the complete collection of the regularized Bayes rules $t_G^*(x; \rho)$ in (3.3) as

$$\mathscr{T}_\rho = \{t_G^*(\cdot; \rho) : G \in \mathscr{G}\}, \tag{4.6}$$

where $\mathscr{G}$ is the family of all distribution functions. The following proposition, proved in the Appendix, provides an entropy bound for (4.6) under the seminorm $\|h\|_{\infty,M} = \sup_{|x|\le M} |h(x)|$.

PROPOSITION 3. *Let* $\widetilde{L}(y) = \sqrt{-\log(2\pi y^2)}$ *be the inverse of* $y = \varphi(x)$ *as in Proposition 1. Then, for all* $0 < \eta \le \rho \le (2\pi e)^{-1/2}$,

$$\begin{aligned} &\log N(\eta^*, \mathscr{T}_\rho, \|\cdot\|_{\infty,M}) \\ &\quad \le \{4(6\widetilde{L}^2(\eta) + 1)(2M/\widetilde{L}(\eta) + 3) + 2\}|\log \eta|, \end{aligned} \tag{4.7}$$

*where* $\eta^* = (\eta/\rho)\{3\widetilde{L}(\eta) + 2\}$.

4.1.4. *An exponential inequality for the loss of regularized Bayes rules.* The last element of our proof is an exponential inequality for the difference between the loss and risk of regularized Bayes rules $t^*_G(\mathbf{X};\rho)$. For each separable rule $t(x)$, the squared loss $\|t(\mathbf{X})-\boldsymbol{\theta}\|^2$ is a sum of independent variables. However, a direct application of the empirical process theory to the loss would yield an oracle inequality of the $n^{-1/2}$ order, which is inadequate for the sharper convergence rates in this paper. Thus, we use the following isoperimetric inequality for the square root of the loss.

PROPOSITION 4. *Suppose $\mathbf{X}\sim N(\boldsymbol{\theta},\mathbf{I}_n)$ under $P_{n,\boldsymbol{\theta}}$. Let $t_G(x;\rho)$ be the regularized Bayes rule as in (3.3), with a deterministic distribution $G$ and $0<\rho\le(2\pi e^3)^{-1/2}$. Let $\widetilde{L}(\rho)=\sqrt{-\log(2\pi\rho^2)}$. Then, for all $x>0$*

$$P_{n,\boldsymbol{\theta}}\{\|t^*_G(\mathbf{X};\rho)-\boldsymbol{\theta}\|\ge E_{n,\boldsymbol{\theta}}\|t^*_G(\mathbf{X};\rho)-\boldsymbol{\theta}\|+x\}\le\exp\biggl(-\frac{x^2}{2\widetilde{L}^4(\rho)}\biggr).$$

PROOF. Let $h(\mathbf{x})=\|t^*_G(\mathbf{x};\rho)-\boldsymbol{\theta}\|$. It follows from Proposition 1 that

$$\begin{aligned}|h(\mathbf{x})-h(\mathbf{y})|&\le\|t^*_G(\mathbf{x};\rho)-t^*_G(\mathbf{y};\rho)\|\\&\le\|\mathbf{x}-\mathbf{y}\|\sup_x|(\partial/\partial x)t^*_G(x;\rho)|\le\widetilde{L}^2(\rho)\|\mathbf{x}-\mathbf{y}\|.\end{aligned}$$

Thus, $h(\mathbf{x})/\widetilde{L}^2(\rho)$ has the unit Lipschitz norm. The conclusion follows from the Gaussian isoperimetric inequality [4]. See page 439 of [34]. □

4.2. *An oracle inequality.* Our oracle inequality for the GMLEB, stated in Theorem 5 below, is a key result of this paper from a mathematical point of view. It builds upon Theorems 3 and 4 and Propositions 2, 3 and 4 (the regularized Bayes rules with misspecified prior, generalized MLE of normal mixtures, representation of the GMLEB, entropy bounds and Gaussian concentration inequality) and leads to adaptive ratio optimality and minimax theorems more general than Theorems 1 and 2.

THEOREM 5. *Let $\mathbf{X}\sim N(\boldsymbol{\theta},\mathbf{I}_n)$ under $P_{n,\boldsymbol{\theta}}$ with a deterministic $\boldsymbol{\theta}\in\mathbb{R}^n$ as in (2.1). Let $L_n(\cdot,\cdot)$ be the average squared loss in (2.2) and $0<p\le\infty$. Let $t^*_{\widehat{G}_n}(\mathbf{X})$ be the GMLEB estimator (2.15) with an approximate generalized MLE $\widehat{G}_n$ satisfying (2.14). Then, there exists a universal constant $M_0$ such that for all $\log n\ge 2/p$,*

$$\begin{aligned}\tilde{r}_{n,\boldsymbol{\theta}}(t^*_{\widehat{G}}(\mathbf{X}))&=\sqrt{E_{n,\boldsymbol{\theta}}L_n(t^*_{\widehat{G}_n}(\mathbf{X}),\boldsymbol{\theta})}-\sqrt{R^*(G_n)}\\&\le M_0\varepsilon_n(\log n)^{3/2},\end{aligned}\tag{4.8}$$

*where $R^*(G_n)$ is the minimum risk of all separable estimators as in (2.8) with $G_n = G_{n,\boldsymbol{\theta}}$ as in (2.4), and $\varepsilon_n = \varepsilon(n, G_n, p)$ is as in (4.2). In particular, for any sequences of constants $M_n \to \infty$ and fixed positive $\alpha$ and $c$,*

$$\varepsilon_n \asymp \begin{cases} n^{-p/(2+2p)}(\log n)^{(2+3p)/(4+4p)}, \\ \qquad \textit{if } \mu_p^w(G_n) = O(1) \textit{ with a fixed } p, \\ n^{-1/2}(\log n)^{3/4}\{M_n^{1/2} \vee (\log n)^{1/4}\}, \\ \qquad \textit{if } G_n([-M_n, M_n]) = 1 \textit{ and } p = \infty, \\ n^{-1/2}(\log n)^{1/(2(2\wedge\alpha))+3/4}, \\ \qquad \textit{if } \int e^{|cu|^\alpha} G_n(du) = O(1) \textit{ and } p \asymp \log n. \end{cases}$$

REMARK 6. In the proof of Theorem 5, applications of Theorems 3 and 4 resulted in the leading term for the upper bound in (4.18), while the contributions of other parts of the proof are of smaller order.

REMARK 7. The $M_0$ in (4.8) is universal for $q_n = (e\sqrt{2\pi}/n^2) \wedge 1$ in (2.14) and depends on $\sup_n |\log q_n|/\log n$ in general.

The consequences of Theorem 5 upon the adaptive ratio optimality and minimaxity of the GMLEB are discussed in the next section. Here is an outline of its proof. The large deviation inequality in Theorem 4 and the representation of the GMLEB in (4.4) imply that

$$\|t^*_{\widehat{G}_n}(\mathbf{X}) - \boldsymbol{\theta}\| \le \|t^*_{\widehat{G}_n}(\mathbf{X}; \rho_n) - \boldsymbol{\theta}\| I_{A_n} + \zeta_{1n}, \qquad \rho_n = \frac{q_n}{e\sqrt{2\pi}n}, \tag{4.9}$$

where $A_n = \{d(f_{\widehat{G}_n}, f_{G_n}) \le x^* \varepsilon_n\}$ and $\zeta_{1n} = \|t^*_{\widehat{G}_n}(\mathbf{X}; \rho_n) - \boldsymbol{\theta}\| I_{A_n^c}$ with $x^* = x_* \vee 1$. By (3.2) and Proposition 1, $|t^*_G(X_i; \rho_n) - \theta_i| \le \widetilde{L}(\rho_n) + |N(0,1)|$, so that Theorem 4 provides an upper bound for $E_{n,\boldsymbol{\theta}}\zeta_{1n}^2$. By the entropy bound in Proposition 3, there exists a finite collection of distributions $\{H_j, j \le N\}$ of manageable size $N$ such that

$$\zeta_{2n} = \left\{\|t^*_{\widehat{G}_n}(\mathbf{X}; \rho_n) - \boldsymbol{\theta}\| I_{A_n} - \max_{j \le N} \|t^*_{H_j}(\mathbf{X}; \rho_n) - \boldsymbol{\theta}\|\right\}_+ \tag{4.10}$$

is small and $d(f_{H_j}, f_{G_n}) \le x^* \varepsilon_n$ for all $j \le N$. Since the regularized Bayes rules $t^*_{H_j}(\mathbf{X}; \rho_n)$ are separable and the collection $\{H_j, j \le N\}$ is of manageable size, the large deviation inequality in Proposition 4 implies that

$$\zeta_{3n} = \max_{j \le N} \{\|t^*_{H_j}(\mathbf{X}; \rho_n) - \boldsymbol{\theta}\| - E_{n,\boldsymbol{\theta}}\|t^*_{H_j}(\mathbf{X}; \rho_n) - \boldsymbol{\theta}\|\}_+ \tag{4.11}$$

is small. Since $d(f_{H_j}, f_{G_n}) \le x^* \varepsilon_n$, Theorem 3 implies that

$$\zeta_{4n} = \max_{j \le N} \sqrt{E_{n,\boldsymbol{\theta}}\|t^*_{H_j}(\mathbf{X}; \rho_n) - \boldsymbol{\theta}\|^2} - \sqrt{nR^*(G_n)} \tag{4.12}$$

is no greater than $O(x^*\varepsilon_n)(\log\rho_n)^{3/2}$, where $R^*(G_n)$ is the general EB benchmark risk in (2.8). Finally, upper bounds for individual pieces $E_{n,\boldsymbol{\theta}}\zeta_{jn}^2$ are put together via

$$\sqrt{E_{n,\boldsymbol{\theta}}\|t^*_{\widehat{G}_n}(\mathbf{X})-\boldsymbol{\theta}\|^2} \le \sqrt{nR^*(G_n)} + \sqrt{E_{n,\boldsymbol{\theta}}\left(\sum_{j=1}^4 |\zeta_{jn}|\right)^2}. \tag{4.13}$$

4.3. *Adaptive ratio optimality and minimaxity.* We discuss here the optimality properties of the GMLEB as consequences of the oracle inequality in Theorem 5.

Theorem 5 immediately implies the adaptive ratio optimality (2.21) of the GMLEB in the classes $\Theta_n^* = \Theta_n^*(M_n)$ for any sequences of constants $M_n \to \infty$, where

$$\Theta_n^*(M) = \left\{\boldsymbol{\theta}\in\mathbb{R}^n : R^*(G_{n,\boldsymbol{\theta}}) \ge M(\log n)^3 \inf_{p\ge 2/\log n} \varepsilon^2(n, G_{n,\boldsymbol{\theta}}, p)\right\} \tag{4.14}$$

with $G_{n,\boldsymbol{\theta}} = G_n$ as in (2.4) and $\varepsilon(n,G,p)$ as in (4.2). This is formally stated in the theorem below.

THEOREM 6. *Let $\mathbf{X}\sim N(\boldsymbol{\theta},\mathbf{I}_n)$ under $P_{n,\boldsymbol{\theta}}$ with a deterministic $\boldsymbol{\theta}\in\mathbb{R}^n$. Let $t^*_{\widehat{G}_n}(\mathbf{X})$ be the GMLEB estimator (2.15) with the approximate MLE $\widehat{G}_n$ in (2.14). Let $R^*(G_{n,\boldsymbol{\theta}})$ be the general EB benchmark in (2.8) with the distribution $G_n = G_{n,\boldsymbol{\theta}}$ in (2.4). Then, for the classes $\Theta_n^*(M)$ in (4.14),*

$$\lim_{(n,M)\to(\infty,\infty)} \sup_{\boldsymbol{\theta}\in\Theta_n^*(M)} \{E_{n,\boldsymbol{\theta}} L_n(t^*_{\widehat{G}_n}(\mathbf{X}),\boldsymbol{\theta})/R^*(G_{n,\boldsymbol{\theta}})\} \le 1. \tag{4.15}$$

REMARK 8. Since the minimum of $\varepsilon(n, G_{n,\boldsymbol{\theta}}, p)$ is taken in (4.14) over $p \ge 2/\log n$ for each $\boldsymbol{\theta}$, the adaptive ratio optimality (4.15) allows smaller $R^*(G_{n,\boldsymbol{\theta}})$ than simply using $\varepsilon(n, G_{n,\boldsymbol{\theta}}, \infty)$ does as in Theorem 1. Thus, Theorem 6 implies Theorem 1.

Another main consequence of the oracle inequality in Theorem 5 is the adaptive minimaxity (2.28) of the GMLEB for a broad range of sequences $\Theta_n \in \mathbb{R}^n$. We have stated our results for regular $\ell_p$ balls in Theorem 2. In the rest of the section, we consider weak $\ell_p$ balls

$$\Theta_{p,C,n}^w = \{\boldsymbol{\theta}\in\mathbb{R}^n : \mu_p^w(G_{n,\boldsymbol{\theta}}) \le C\}, \tag{4.16}$$

where $G_{n,\boldsymbol{\theta}}$ is the empirical distribution of the components of $\boldsymbol{\theta}$ and the functional $\mu_p^w(G)$ is the weak moment in (4.1). Alternatively,

$$\Theta_{p,C,n}^w = \left\{\boldsymbol{\theta}\in\mathbb{R}^n : \max_{1\le i\le n} |\theta_i|^p \sum_{j=1}^n I\{|\theta_j|\ge|\theta_i|\}/n \le C^p\right\}.$$

THEOREM 7. *Let $\mathbf{X} \sim N(\boldsymbol{\theta}, \mathbf{I}_n)$ under $P_{n,\boldsymbol{\theta}}$ with a deterministic $\boldsymbol{\theta} \in \mathbb{R}^n$. Let $L_n(\cdot,\cdot)$ be the average squared loss (2.2) and $\mathscr{R}_n(\Theta)$ be the minimax risk (2.27). Then, for all approximate solutions $\widehat{G}_n$ satisfying (2.14), the GMLEB $\widehat{\boldsymbol{\theta}} = t^*_{\widehat{G}_n}(\mathbf{X})$ is adaptive minimax (2.28) in the weak $\ell_p$ balls $\Theta_n = \Theta^w_{p,C_n,n}$ in (4.16), provided that the radii $C_n$ are within the range (2.30).*

Here is our argument. The weak $L_p$ ball that matches (4.16) is

$$\mathscr{G}^w_{p,C} = \{G : \mu^w_p(G) \le C\}.$$

Let $J^w_{p,C}(\lambda) = -\int_0^\infty (t^2 \wedge \lambda^2) d\{1 \wedge (C/t)^p\}$, which is approximately the Bayes risk of the soft-threshold estimator for the stochastically largest Pareto prior in $\mathscr{G}_{p,C}$. Let $\lambda_{p,C} = \sqrt{1 \vee \{2\log(1/C^{p\wedge 2})\}}$. Johnstone [23] proved that

$$\lim_{n\to\infty} \frac{\mathscr{R}_n(\Theta^w_{p,C_n,n})}{\mathscr{R}(\mathscr{G}^w_{p,C_n})} = 1 \tag{4.17}$$

for $p > 2$ with $C_n \to C+ \ge 0$ and for $p \le 2$ with $nC_n^p/(\log n)^{1+6/p} \to \infty$, and that $\mathscr{R}(\mathscr{G}^w_{p,C_n})/J^w_{p,C_n}(\lambda_{p,C_n}) \to 1$ as $C_n^{p\wedge 2} \to 0$. Abramovich et al. [1] proved $\mathscr{R}_n(\Theta^w_{p,C_n,n})/J^w_{p,C_n}(\lambda_{p,C_n}) \to 1$ for $p < 2$ and $(\log n)^5/n \le C_n^p \le n^{-\kappa}$ for all $\kappa > 0$. The combination of their results implies (4.17) for $p \le 2$ and $C_n^p \ge (\log n)^5/n$. Therefore, (4.17) holds under (2.30) due to $pk_1(p) = p/2 + 4 + 3/p > 5$ for $p < 2$. As in Section 2.7,

$$\sup_{\boldsymbol{\theta}\in\Theta^w_{p,C_n,n}} \widetilde{r}_{n,\boldsymbol{\theta}}(\widehat{t}_n) = o(1)\sqrt{J_{p,C_n}} \tag{4.18}$$

as in (2.37), due to $J_{p,C_n} \asymp \mathscr{R}(\mathscr{G}_{p,C_n}) \le \mathscr{R}(\mathscr{G}^w_{p,C_n})$.

**5. More simulation results.** In addition to the simulation results reported in Section 2.4, we conducted more experiments to explore a larger sample size, sparse unknown means without exact zero, and i.i.d. unknown means from normal priors. The results for the nine statistical procedures and the oracle rule $t^*_{G_n}(\mathbf{X})$ for the general EB are reported in Tables 2–4, in the same format as Table 1. Each entry is based on an average of 100 replications. In each column, boldface entries indicate top three performers other than the hybrid estimator or the oracle. Two columns with $\mu = 4$ are dropped to fit the tables in.

In Table 2 we report simulation results for $n = 4000$. Compared with Table 1, F-GEB replaces EBThresh as a distant third top performer in the moderately sparse case of $\#\{i : \theta_i = \mu\} = 200$, and almost the same sets of estimators prevail as top performers in other columns. Since the collections of $G_n$ are identical in Tables 1 and 2, the average squared loss $\|\widehat{\boldsymbol{\theta}} - \boldsymbol{\theta}\|^2/n$ should decrease in $n$ to indicate convergence to the oracle risks for each estimator in each model, but this is not the case in entries in italics.

TABLE 2
*Average of* $\|\widehat{\boldsymbol{\theta}} - \boldsymbol{\theta}\|^2 : n = 4000,\ \theta_i \in \{0, \mu\},\ \#\{i : \theta_i = \mu\} =$ *20, 200 or 2000*

| **# nonzero** | **20** | | | | **200** | | | **2000** | | |
|---|---|---|---|---|---|---|---|---|---|---|
| $\boldsymbol{\mu}$ | **3** | **4** | **5** | **7** | **3** | **5** | **7** | **3** | **5** | **7** |
| James–Stein | 175 | 298 | 446 | 790 | 1243 | *2229* | 2846 | 3261 | 3689 | *3829* |
| EBThresh | **145** | **120** | **63** | *377* | *861* | 404 | 290 | 3411 | 3118 | *2621* |
| SURE | *174* | *270* | *329* | *355* | *1725* | 827 | 827 | 3296 | 3317 | *3317* |
| FDR (0.01) | *175* | 202 | 103 | **26** | *1569* | 506 | *231* | 10,230 | 2607 | *2090* |
| FDR (0.1) | 161 | 138 | 70 | 48 | *1121* | 450 | *409* | 4578 | 2597 | *2563* |
| GMLEB | **141** | **115** | **68** | **30** | **624** | **215** | **43** | **1808** | **489** | **62** |
| S-GMLEB | **116** | **92** | **45** | **10** | **597** | **193** | **23** | **1791** | **479** | **53** |
| F-GEB | 243 | 231 | 166 | 156 | **739** | **353** | **229** | **1907** | **641** | **253** |
| HF-GEB | 145 | 120 | 63 | *377* | 694 | 286 | 159 | 1868 | 576 | 171 |
| Oracle | 110 | 84 | 40 | 3 | 587 | 186 | 16 | 1771 | 460 | 36 |

In Table 3, we report simulation results for sparse mean vectors without exact zero. It turns out that adding uniform $[-0.2, 0.2]$ perturbations to $\theta_i$ does not change the results much, compared with Table 1.

In Table 4, we report simulation results for i.i.d. $\theta_i \sim N(\mu, \sigma^2)$. This is the parametric model in which the (oracle) Bayes estimators are linear. Indeed, the James–Stein estimator is the top performer throughout all the columns and tracks the oracle risk extremely well, while the GMLEB is not so far behind. It is interesting that for $\sigma^2 = 40$, the EBThresh and SURE outperform GMLEB as they approximate the naive $\widehat{\boldsymbol{\theta}} = \mathbf{X}$ with diminishing threshold levels. Another interesting phenomenon is the disappearance of the advantage of the S-GMLEB over the GMLEB, as the unknowns are no longer sparse.

TABLE 3
*Average of* $\|\widehat{\boldsymbol{\theta}} - \boldsymbol{\theta}\|^2$*:* $n = 1000,\ \theta_i = \mu_i + \mathrm{unif}[-0.2, 0.2],\ \mu_i \in \{0, \mu\},\ \#\{i : \mu_i = \mu\} =$*5, 50 or 500*

| $\boldsymbol{\#\{\mu_i \neq 0\}}$ | **5** | | | | **50** | | | **500** | | |
|---|---|---|---|---|---|---|---|---|---|---|
| $\boldsymbol{\mu}$ | **3** | **4** | **5** | **7** | **3** | **5** | **7** | **3** | **5** | **7** |
| James–Stein | 57 | 87 | 124 | 207 | 316 | 559 | 713 | 817 | 932 | 971 |
| EBThresh | **48** | **44** | **31** | **23** | **226** | **115** | 87 | 855 | 797 | 677 |
| SURE | 55 | 75 | 84 | 89 | 426 | 221 | 220 | 830 | 845 | 848 |
| FDR (0.01) | 56 | 62 | 37 | **20** | 395 | 137 | **72** | 2555 | 676 | 541 |
| FDR (0.1) | 53 | 49 | 34 | 27 | 289 | 130 | 116 | 1152 | 666 | 664 |
| GMLEB | **49** | **45** | **32** | **23** | **170** | **70** | **27** | **466** | **146** | **32** |
| S-GMLEB | **45** | **41** | **29** | **19** | **164** | **67** | **24** | **462** | **145** | **31** |
| F-GEB | 115 | 108 | 105 | 91 | 238 | 155 | 118 | **534** | **244** | **145** |
| HF-GEB | 48 | 44 | 31 | 23 | 210 | 113 | 85 | 509 | 203 | 101 |
| Oracle | 39 | 35 | 23 | 14 | 158 | 61 | 18 | 454 | 135 | 22 |

TABLE 4
*Average of* $\|\widehat{\boldsymbol{\theta}} - \boldsymbol{\theta}\|^2$: $n = 1000$, *i.i.d.* $\theta_j \sim N(\mu, \sigma^2)$

| $\boldsymbol{\sigma^2}$ | **0.1** | | | | **2** | | | **40** | | |
|---|---|---|---|---|---|---|---|---|---|---|
| $\boldsymbol{\mu}$ | **3** | **4** | **5** | **7** | **3** | **5** | **7** | **3** | **5** | **7** |
| James–Stein | **92** | **92** | **92** | **93** | **665** | **670** | **665** | **970** | **982** | **975** |
| EBThresh | 1081 | 1058 | 1035 | 1020 | 1013 | 1032 | 1014 | **983** | **998** | **997** |
| SURE | 1006 | 1505 | 3622 | 13,146 | 988 | 1033 | 3514 | **983** | **998** | **996** |
| FDR (0.01) | 3972 | 2049 | 1169 | 999 | 2789 | 1599 | 1050 | 1661 | 1566 | 1427 |
| FDR (0.1) | 1555 | 1093 | 1002 | 998 | 1455 | 1096 | 999 | 1184 | 1161 | 1117 |
| GMLEB | **94** | **94** | **95** | **95** | **675** | **678** | **673** | 1001 | 1015 | 1009 |
| S-GMLEB | **97** | **98** | **99** | **98** | **678** | **681** | **675** | 1002 | 1015 | 1009 |
| F-GEB | 171 | 171 | 175 | 171 | 735 | 743 | 736 | 1107 | 1130 | 1122 |
| HF-GEB | 138 | 139 | 143 | 142 | 721 | 726 | 720 | 1067 | 1088 | 1079 |
| Oracle | 91 | 90 | 91 | 90 | 665 | 669 | 664 | 970 | 981 | 975 |

**6. Discussion.** In this section, we discuss general EB with kernel estimates of the oracle Bayes rule, sure computation of an approximate generalized MLE and a number of additional issues.

6.1. *Kernel methods.* General EB estimators of the mean vector $\boldsymbol{\theta}$ can be directly derived from the formula (3.2) using the kernel method

$$\widehat{\boldsymbol{\theta}} = \widehat{t}_n(\mathbf{X}), \qquad \widehat{t}_n(x) = x + \frac{\widehat{f}_n'(x)}{\widehat{f}_n(x) \vee \rho_n}, \tag{6.1}$$
$$\widehat{f}_n(x) = \sum_{i=1}^{n} \frac{K(a_n(x - X_i))}{n a_n}.$$

This was done in [36] with the Fourier kernel $K(x) = (\sin x)/(\pi x)$ and $\sqrt{2 \log n} \le a_n \asymp \sqrt{\log n}$. The main rationale for using the Fourier kernel is the near optimal convergence rate of $\widehat{f}_n - f_{G_n} = O(\sqrt{(\log n)/n})$ and $\widehat{f}_n' - f_{G_n}' = O((\log n)/\sqrt{n})$, uniformly in $\boldsymbol{\theta}$. However, since the relationship between $\widehat{f}_n'(x)$ and $\widehat{f}_n(x)$ is not as trackable as in the case of generalized MLE $f_{\widehat{G}_n}$, a much higher regularization level $\rho_n \asymp \sqrt{(\log n)/n}$ in (6.1) were used [36, 38] to justify the theoretical results. This could be an explanation for the poor performance of the Fourier general EB estimator for very sparse $\boldsymbol{\theta}$ in our simulations. From this point of view, the GMLEB is much more appealing since its estimating function retains all analytic properties of the Bayes rule. Consequently, the GMLEB requires no regularization for the adaptive ratio optimality and adaptive minimaxity in our theorems.

Brown and Greenshtein [6] have studied (6.1) with the normal kernel $K(x) = \varphi(x)$ and possibly different bandwidth $1/a_n$, and have proved the adaptive ratio optimality (2.21) of their estimator when $\|\boldsymbol{\theta}\|_\infty$ and $R^*(G_{n,\boldsymbol{\theta}})$ have certain different polynomial orders. The estimating function $\widehat{t}_n(x)$ with the normal kernel, compared with the Fourier kernel, behaves more like the regularized Bayes rule (3.3) analytically with the positivity of $\widehat{f}_n(x)$ and more trackable relationship between $\widehat{f}'_n(x)$ and $\widehat{f}_n(x)$. Still, it is unclear without some basic properties of the Bayes rule in Proposition 1 and Theorem 3, it is unclear if the kernel methods of the form (6.1) would possess as strong theoretical properties as in Theorems 1, 2, 5, 6 and 7 or perform as well as the GMLEB for moderate samples in simulations [6].

6.2. *Sure computation of an approximate general MLE.* We present a conservative data-driven criterion to guarantee (2.14) with the EM-algorithm. This provides a definitive way of computing the map from $\{X_i\}$ to $\widehat{G}$ in (2.14) and then to the GMLEB via (2.18).

Set $u_1 = \min_{1\le i\le n} X_i$, $u_m = \max_{1\le i\le n} X_i$, and

$$\varepsilon = (u_m - u_1)/(m-1), \qquad u_j = u_{j-1} + \varepsilon. \tag{6.2}$$

PROPOSITION 5. *Suppose $\varepsilon^2\{(u_m - u_1)^2/4 + 1/8\} \le 1/n$ with a sufficiently large $m$ in (6.2). Let $w_j^{(0)} > 0$ $\forall j \le m$ with $\sum_{j=1}^m \widehat{w}_j^{(0)} = 1$. Suppose that the EM-algorithm (2.19) is stopped at or beyond an iteration $k > 0$ with*

$$\max_{1\le j\le m} \log(\widehat{w}_j^{(k)}/\widehat{w}_j^{(k-1)}) \le \frac{1}{n}\log\left(\frac{1}{eq_n}\right). \tag{6.3}$$

*Then, (2.14) holds for $\widehat{G}_n = \sum_{j=1}^m \widehat{w}_j^{(k)}\delta_{u_j}$.*

Heuristically, smaller $m$ provides larger $\min_j \widehat{w}_j^{(k)}$ and faster convergence of the EM-algorithm, so that the "best choice" of $m$ is

$$m - 2 < (u_m - u_1)\sqrt{n\{(u_m - u_1)^2/4 + 1/8\}} \le m - 1.$$

For $\max_i |X_i| \asymp \sqrt{\log n}$, this ensures the first condition of Proposition 5 with $m \asymp (\log n)\sqrt{n}$ and $\varepsilon \asymp (n\log n)^{-1/2}$. Proposition 5 is proved via the smoothness of the normal density and Cover's upper bound [10, 35] for the maximum likelihood in finite mixture models.

6.3. *Additional remarks.* A crucial element for the theoretical results in this paper is the oracle inequality for the regularized Bayes estimator with misspecified prior, as stated in Theorem 3. However, we do not believe that mathematical induction is sharp in the argument with higher and higher order of differentiation in the proof of Lemma 1. Consequently, the power $\kappa_1$ in Theorems 2 and 7 is larger than its counterpart more directly established for threshold estimators [1, 24]. Still, the GMLEB performs as well as any threshold estimators in our simulations for the most sparse mean vectors. As expected, the gain of the GMLEB is huge against the James–Stein estimator for sparse means and against threshold estimators for dense means.

It is interesting to observe in Tables 1–3 that the simulated $\ell_2$ risk for the GMLEB sometimes dips well below the benchmark $\sum_{i=1}^{n}\theta_i^2 \wedge 1 = \#\{i \le n : \theta_i \ne 0\}$ for the oracle hard threshold rule $\widehat{\theta}_i = X_i I\{|\theta_i| \le 1\}$ [18], while the simulated $\ell_2$ risk for threshold estimators is always above that benchmark.

An important consequence of our results is the adaptive minimaxity and other optimality properties of the GMLEB approach to nonparametric regression under suitable smoothness conditions. For example, applications of the GMLEB estimator to the observed wavelet coefficients at individual resolution levels yield adaptive exact minimaxity in all Basov balls as in [38].

The adaptive minimaxity (2.28) in Theorems 2 and 7 is uniform in the radii $C$ for fixed shape $p$. A minimax theory for (weak) $\ell_p$ balls uniform in $(p, C)$ can be developed by careful combination and improvement of the proofs in [12, 23, 38]. Since the oracle inequality (4.8) is uniform in $p$, uniform adaptive minimaxity in both $p$ and $C$ is in principle attainable for the GMLEB.

The theoretical results in this paper are all stated for deterministic $\boldsymbol{\theta} = (\theta_1, \dots, \theta_n)$. By either mild modifications of the proofs here or conditioning on the unknowns, analogues versions of our theorems can be established for the estimation of i.i.d. means $\{\xi_i\}$ in the EB model (2.25). Other possible directions of extension of the results in this paper are the cases of $X_i \sim N(\theta_i, \sigma_n^2)$ via scale change, with known $\sigma_n^2$ or an independent consistent estimate of $\sigma_n^2$, and $X_i \sim N(\theta_i, \sigma_i^2)$ with known $\sigma_i^2$.

## APPENDIX

Here we prove Proposition 1, Lemma 1, Proposition 3, Theorems 5, 2 and 7, and then Proposition 5. We need one more lemma for the proof of Proposition 1. Throughout this appendix, $\lfloor x \rfloor$ denotes the greatest integer lower bound of $x$, and $\lceil x \rceil$ denotes the smallest integer upper bound of $x$.

LEMMA A.1. *Let $f_G(x)$ be as in (2.13) and $\widetilde{L}(y)$ as in Proposition 1. Then,*

$$\left(\frac{f_G'(x)}{f_G(x)}\right)^2 \le \frac{f_G''(x)}{f_G(x)} + 1 \le \widetilde{L}^2(f_G(x)) = \log\left(\frac{1}{2\pi f_G^2(x)}\right) \qquad \forall x. \tag{A.1}$$

PROOF. Since $Y|\xi \sim N(\xi,1)$ and $\xi \sim G$ under $P_G$, by (3.2)

$$\frac{f_G'(x)}{f_G(x)} = E_G[\xi - Y|Y=x],$$

$$\frac{f_G''(x)}{f_G(x)} + 1 = E_G[(\xi - Y)^2|Y=x].$$

This gives the first inequality of (A.1). The second inequality of (A.1) follows from Jensen's inequality: for $h(x) = e^{x/2}$

$$h\left(\frac{f_G''(x)}{f_G(x)} + 1\right) \le E_G[h((\xi - Y)^2)|Y=x] = \frac{1}{\sqrt{2\pi} f_G(x)}.$$

This completes the proof. □

PROOF OF PROPOSITION 1. Since $f_G(x) = \int \varphi(x-u)G(du) \ge 0$, the value of (3.3) is always between $t_G^*(x)$ and $x$. By Lemma A.1

$$|x - t_G^*(x;\rho)| \le \frac{f_G(x)}{f_G(x) \vee \rho} \widetilde{L}(f_G(x)) \le \widetilde{L}(\rho)$$

for $\rho \le (2\pi e)^{-1/2}$, since $\widetilde{L}(y)$ is decreasing in $y^2$ and $y^2\widetilde{L}^2(y)$ is increasing in $y^2 \le 1/(2\pi e)$. Similarly, the second line of (3.4) follows from Lemma A.1 and

$$\frac{\partial t_G^*(x;\rho)}{\partial x} = \begin{cases} 1 + f_G''(x)/f_G(x) - \{f_G'(x)/f_G(x)\}^2, & f_G(x) > \rho, \\ 1 + f_G''(x)/\rho, & f_G(x) < \rho. \end{cases}$$

Note that $\widetilde{L}(f_G(x)) \le \widetilde{L}(\rho)$ for $f_G(x) \ge \rho$, and for $f_G(x) < \rho \le (2\pi e^3)^{-1/2}$

$$0 \le 1 - \frac{f_G(x)}{\rho} \le 1 + \frac{f_G''(x)}{\rho} \le 1 + \frac{f_G(x)}{\rho}(\widetilde{L}^2(f_G(x)) - 1) \le \widetilde{L}^2(\rho)$$

due to the monotonicity of $y\{\widetilde{L}^2(y) - 1\}$ in $0 \le y \le (2\pi e^3)^{-1/2}$. □

PROOF OF LEMMA 1. Let $D = d/dx$. We first prove that for all integers $k \ge 0$ and $a \ge \sqrt{2k-1}$,

$$\int \{D^k(f_G - f_{G_0})\}^2 \, dx \le \frac{4a^{2k}}{\sqrt{2\pi}} d^2(f_G, f_{G_0}) + \frac{4a^{2k-1}}{\pi} e^{-a^2}. \tag{A.2}$$

Let $h^*(u) = \int e^{iux} h(x)\, dx$ for all integrable $h$. Since $|f_G^*(u)| \le \varphi^*(u) = e^{-u^2/2}$, it follows from the Plancherel identity that

$$\begin{aligned} \int \{D^k(f_G - f_{G_0})\}^2 \, dx &= \frac{1}{2\pi} \int u^{2k} |f_G^*(u) - f_{G_0}^*(u)|^2 \, du \\ &\le \frac{a^{2k}}{2\pi} \int |f_G^*(u) - f_{G_0}^*(u)|^2 \, du + \frac{4}{2\pi} \int_{|u|>a} u^{2k} e^{-u^2} \, du \\ &= a^{2k} \int |f_G - f_{G_0}|^2 \, dx + \frac{4}{\pi} c_k, \end{aligned}$$

where $c_k = \int_{u>a} u^{2k} e^{-u^2}\,du$. Since $(k-1/2) \le a^2/2$, integrating by parts yields

$$\begin{aligned} c_k &= 2^{-1} a^{2k-1} e^{-a^2} + \{(k-1/2)/a^2\} a^2 c_{k-1} \\ &\le 2^{-1} a^{2k-1} e^{-a^2} (1 + 1/2 + \cdots + 1/2^{k-1}) + 2^{-k} a^{2k} c_0 \\ &\le a^{2k-1} e^{-a^2} \end{aligned}$$

due to $c_0 \le a^{-1} \int_{u>a} u e^{-u^2}\,du = e^{-a^2}/(2a)$. Since $f_G(x) \le 1/\sqrt{2\pi}$,

$$\int |f_G - f_{G_0}|^2\,dx \le \|\sqrt{f_G} + \sqrt{f_{G_0}}\|_\infty^2 d^2(f_G, f_{G_0}) \le \frac{4}{\sqrt{2\pi}} d^2(f_G, f_{G_0}).$$

The combination of the above inequalities yields (A.2).

Define $w_* = 1/(f_G \vee \rho + f_{G_0} \vee \rho)$ and $\Delta_k = (\int \{D^k(f_G - f_{G_0})\}^2 w_*)^{1/2}$. Integrating by parts, we find

$$\Delta_k^2 = -\int \{D^{k-1}(f_G - f_{G_0})\}\{D^{k+1}(f_G - f_{G_0}) w_* + (D^k(f_G - f_{G_0}))(D w_*)\}.$$

Since $|(Dw_*)(x)| \le 2\widetilde{L}(\rho) w_*(x)$ by Proposition 1, Cauchy–Schwarz gives

$$\Delta_k^2 \le \Delta_{k-1}\Delta_{k+1} + 2\widetilde{L}(\rho)\Delta_{k-1}\Delta_k.$$

Let $k_0$ be a nonnegative integer satisfying $k_0 \le \widetilde{L}^2(\rho)/2 < k_0 + 1$. Define $k^* = \min\{k : \Delta_{k+1} \le k_0 2\widetilde{L}(\rho)\Delta_k\}$. For $k < k^*$, we have $\Delta_k^2 \le (1 + 1/k_0)\Delta_{k-1}\Delta_{k+1}$, so that for $k^* \le k_0$,

$$\frac{\Delta_1}{\Delta_0} \le \left(1 + \frac{1}{k_0}\right)\frac{\Delta_2}{\Delta_1} \le \left(1 + \frac{1}{k_0}\right)^{k^*} \frac{\Delta_{k^*+1}}{\Delta_{k^*}} \le e k_0 2\widetilde{L}(\rho) \le e\widetilde{L}^3(\rho).$$

Since $(f_G^{1/2} + f_{G_0}^{1/2})^2 w_* \le 2$, we have $\Delta_0^2 \le 2d^2(f_G, f_{G_0})$. Thus, for $k^* \le k_0$

$$\Delta_1 \le e\widetilde{L}^3(\rho)\sqrt{2}\, d(f_G, f_{G_0}). \tag{A.3}$$

For $k_0 < k^*$, $\Delta_1/\Delta_0 \le (1 + 1/k_0)^k \Delta_{k+1}/\Delta_k$ for all $k \le k_0$, so that

$$\begin{aligned} \frac{\Delta_1}{\Delta_0} &\le \left[\prod_{k=0}^{k_0} \{(1 + 1/k_0)^k \Delta_{k+1}/\Delta_k\}\right]^{1/(k_0+1)} \\ &= (1 + 1/k_0)^{k_0/2} \{\Delta_{k_0+1}/\Delta_0\}^{1/(k_0+1)}. \end{aligned} \tag{A.4}$$

To bound $\Delta_{k_0+1}$ by (A.2), we pick the constant $a > 0$ with the $a^2$ in (3.9), so that $a^2 \ge 2(k_0 + 1/2)$ and $e^{-a^2} \le d^2(f_G, f_{G_0})$. Since $w_* \le 1/(2\rho)$, an application of (A.2) with this $a$ gives

$$\begin{aligned} \Delta_{k_0+1}^2 &\le \frac{1}{2\rho} \int \{D^{k_0+1}(f_G - f_{G_0})\}^2 \\ &\le \frac{2a^{2(k_0+1)}}{\rho\sqrt{2\pi}} d^2(f_G, f_{G_0})(1 + a^{-1}\sqrt{2/\pi}). \end{aligned}$$

Since $\Delta_0^2 \le 2d^2(f_G, f_{G_0})$, inserting the above inequality into (A.4) yields

$$\begin{aligned}\Delta_1 &\le (1+1/k_0)^{k_0/2}\Delta_0^{k_0/(k_0+1)}\Delta_{k_0+1}^{1/(k_0+1)}\\ &\le (1+1/k_0)^{k_0/2}\sqrt{2}d(f_G, f_{G_0})a\left(\frac{1+\sqrt{2/\pi}}{\rho\sqrt{2\pi}}\right)^{1/(2k_0+2)} \qquad \text{(A.5)}\\ &\le \sqrt{e2}d(f_G, f_{G_0})a\sqrt{2}(2\pi\rho^2)^{-1/(4k_0+4)}.\end{aligned}$$

Since $|\log(2\pi\rho^2)| = \widetilde{L}^2(\rho) < 2k_0+2$, (3.9) follows from (A.3) and (A.5). $\square$

Proof of Proposition 3. We provide a dense version of the proof since it is similar to the entropy calculations in [20, 21, 39].

It follows from (3.3), (3.4) and Lemma A.1 that

$$|t_G^*(x;\rho) - t_H^*(x;\rho)| \le \frac{1}{\rho}|f_G'(x) - f_H'(x)| + \frac{\widetilde{L}(\rho)}{\rho}|f_G(x) - f_H(x)|, \qquad \text{(A.6)}$$

so that we need to control the norm of both $f_G$ and $f_G'$.

Let $a = \widetilde{L}(\eta)$, $j^* = \lceil 2M/a+2\rceil$ and $k^* = \lfloor 6a^2 \rfloor$. Define semiclosed intervals

$$I_j = (-M+(j-2)a, (-M+(j-1)a)\wedge(M+a)], \qquad j=1,\dots,j^*,$$

to form a partition of $(-M-a, M+a]$. It follows from the Carathéodory's theorem [9] that for each distribution function $G$ there exists a discrete distribution function $G_m$ with support $[-M-a, M+a]$ and no more than $m = (2k^*+2)j^*+1$ support points such that

$$\int_{I_j} u^k G(du) = \int_{I_j} u^k G_m(du), \qquad k = 0,1,\dots,2k^*+1,\ j=1,\dots,j^*. \qquad \text{(A.7)}$$

Since the Taylor expansion of $e^{-t^2/2}$ has alternating signs, for $t^2/2 \le k^*+2$

$$0 \le \operatorname{Rem}(t) = (-1)^{k^*+1}\left\{\varphi(t) - \sum_{k=0}^{k^*}\frac{(-t^2/2)^k}{k!\sqrt{2\pi}}\right\} \le \frac{(t^2/2)^{k^*+1}}{(k^*+1)!\sqrt{2\pi}}.$$

Thus, since $k^*+1 \ge 6a^2$, for $x \in I_j \cap [-M, M]$, the Stirling formula yields

$$\begin{aligned}&|f_G'(x) - f_{G_m}'(x)|\\ &\quad\le \left|\int_{(I_{j-1}\cup I_j\cup I_{j+1})^c}(x-u)\varphi(x-u)\{G(du)-G_m(du)\}\right|\\ &\qquad + \left|\int_{I_{j-1}\cup I_j\cup I_{j+1}}(x-u)\operatorname{Rem}(x-u)\{G(du)-G_m(du)\}\right| \qquad \text{(A.8)}\\ &\quad\le \max_{t\ge a} t\varphi(t) + \frac{4a\{(2a)^2/2\}^{k^*+1}}{\sqrt{2\pi}(k^*+1)!} \le a\eta + \frac{4a(e/3)^{k^*+1}}{2\pi(k^*+1)^{1/2}}\end{aligned}$$

due to $a \ge 1$. Similarly, for $|x| \le M$

$$|f_G(x) - f_{G_m}(x)| \le \eta + \frac{(e/3)^{k^*+1}}{2\pi(k^*+1)^{1/2}}. \tag{A.9}$$

Furthermore, since $(e/3)^6 \le e^{-1/2}$ and $k^*+1 \ge 6a^2 \ge 6$, we have $(e/3)^{k^*+1} \le e^{-a^2/2} = \sqrt{2\pi}\eta$, so that by (A.6), (A.8) and (A.9)

$$\begin{aligned}
&\|t^*_G(\cdot;\rho) - t^*_{G_m}(\cdot;\rho)\|_{\infty,M} \\
&\quad \le \rho^{-1}\biggl(a\eta + \frac{4ae^{-a^2/2}}{2\pi\sqrt{6a^2}}\biggr) + \rho^{-1}\widetilde{L}(\rho)\biggl(\eta + \frac{e^{-a^2/2}}{2\pi\sqrt{6a^2}}\biggr) \qquad \text{(A.10)}\\
&\quad \le \rho^{-1}\eta(2\widetilde{L}(\eta) + 5/\sqrt{12\pi}).
\end{aligned}$$

Let $\xi \sim G_m$, $\xi_\eta = \eta\,\mathrm{sgn}(\xi)\lfloor|\xi|/\eta\rfloor$ and $G_{m,\eta} \sim \xi_\eta$. Since $|\xi - \xi_\eta| \le \eta$,

$$\|f_{G_m} - f_{G_{m,\eta}}\|_\infty \le C_1^*\eta, \qquad \|f'_{G_m} - f'_{G_{m,\eta}}\|_\infty \le C_2^*\eta,$$

where $C_1^* = \sup_x|\varphi'(x)| = (2e\pi)^{-1/2}$ and $C_2^* = \sup_x|\varphi''(x)| = \sqrt{2/\pi}e^{-3/2}$. This and (A.6) imply

$$\|t^*_{G_m}(\cdot;\rho) - t^*_{G_{m,\eta}}(\cdot;\rho)\|_\infty \le \frac{\eta}{\rho}\{C_2^* + C_1^*\widetilde{L}(\rho)\}. \tag{A.11}$$

Moreover, $G_{m,\eta}$ has at most $m$ support points.

Let $\mathscr{P}^m$ be the set of all vectors $\mathbf{w} = (w_1, \dots, w_m)$ satisfying $w_j \ge 0$ and $\sum_{j=1}^m w_j = 1$. Let $\mathscr{P}^{m,\eta}$ be an $\eta$-net of $N(\eta, \mathscr{P}^m, \|\cdot\|_1)$ elements in $\mathscr{P}^m$:

$$\inf_{\mathbf{w}^{m,\eta}\in\mathscr{P}^{m,\eta}} \|\mathbf{w} - \mathbf{w}^{m,\eta}\|_1 \le \eta \qquad \forall \mathbf{w} \in \mathscr{P}^m.$$

Let $\{u_j, j = 1, \dots, m\}$ be the support of $G_{m,\eta}$ and $\mathbf{w}^{m,\eta}$ be a vector in $\mathscr{P}^{m,\eta}$ with $\sum_{j=1}^m |G_{m,\eta}(\{u_j\}) - w_j^{m,\eta}| \le \eta$. Set $\widetilde{G}_{m,\eta} = \sum_{j=1}^m w_j^{m,\eta}\delta_{u_j}$. Then,

$$\|f_{G_{m,\eta}} - f_{\widetilde{G}_{m,\eta}}\|_\infty \le C_0^*\eta, \qquad \|f'_{G_{m,\eta}} - f'_{\widetilde{G}_{m,\eta}}\|_\infty \le C_1^*\eta,$$

where $C_0^* = \varphi(0) = 1/\sqrt{2\pi}$. This and (A.6) imply

$$\|t^*_{G_{m,\eta}}(\cdot;\rho) - t^*_{\widetilde{G}_{m,\eta}}(\cdot;\rho)\|_\infty \le \frac{\eta}{\rho}\{C_1^* + C_0^*\widetilde{L}(\rho)\}. \tag{A.12}$$

The support of $G_{m,\eta}$ and $\widetilde{G}_{m,\eta}$ is $\Omega_{\eta,M} = \{0, \pm\eta, \pm 2\eta, \dots\} \cap [-M-a, M+a]$.

Summing (A.10), (A.11) and (A.12) together, we find

$$\begin{aligned}
&\|t^*_G(\cdot;\rho) - t^*_{\widetilde{G}_{m,\eta}}(\cdot;\rho)\|_{\infty,M} \\
&\quad \le (\eta/\rho)[\{2 + C_1^* + C_0^*\}\widetilde{L}(\eta) + 5/\sqrt{12\pi} + C_2^* + C_1^*] \\
&\quad \le (\eta/\rho)\{2.65\widetilde{L}(\eta) + 1.24\} \le \eta^*.
\end{aligned}$$

Counting the number of ways to realize $\{u_j\}$ and $\mathbf{w}^{m,\eta}$, we find

$$N(\eta^*, \mathscr{T}_\rho, \|\cdot\|_{\infty,M}) \le \binom{|\Omega_{\eta,M}|}{m} N(\eta, \mathscr{P}^m, \|\cdot\|_1), \tag{A.13}$$

with $m = (2k^*+2)j^*+1$, $|\Omega_{\eta,M}| = 1+2\lfloor (M+a)/\eta\rfloor$, $a = \widetilde{L}(\eta)$, $j^* = \lceil 2M/a + 2\rceil$ and $k^* = \lfloor 6a^2\rfloor$.

Since $\mathscr{P}^m$ is in the $\ell_1$ unit-sphere of $\mathbb{R}^m$, $N(\eta, \mathscr{P}^m, \|\cdot\|_1)$ is no greater than the maximum number of disjoint $\mathrm{Ball}(\mathbf{v}_j, \eta/2, \|\cdot\|_1)$ with centers $\mathbf{v}_j$ in the unit sphere. Since all these balls are inside the $(1+\eta/2)$ $\ell_1$-ball, volume comparison yields $N(\eta, \mathscr{P}^m, \|\cdot\|_1) \le (2/\eta+1)^m$. With another application of the Stirling formula, this and (A.13) yield

$$\begin{aligned} N(\eta^*, \mathscr{T}_\rho, \|\cdot\|_{\infty,M}) &\le (2/\eta+1)^m |\Omega_{\eta,M}|^m / m! \\ &\le \{(1+2/\eta)(1+2(M+a)/\eta)\}^m \{(m+1)^{m+1/2} e^{-m-1}\sqrt{2\pi}\}^{-1} \\ &\le [(\eta+2)(\eta+2(M+a))e/(m+1)]^m \eta^{-2m} e\{2\pi(m+1)\}^{-1/2}. \end{aligned} \tag{A.14}$$

Since $m-1 \ge 12a^2(2M/a+2) = 24a(M+a)$ and $a \ge 1 \ge 1/2 \ge \eta$,

$$(\eta+2)(\eta+2(M+a))e \le 8\{1/2+2(M+a)\} \le m+1.$$

Hence, (A.14) is bounded by $\eta^{-2m}$ with $m \le 2(6a^2+1)(2M/a+3)+1$. □

Proof of Theorem 5. Throughout the proof, we use $M_0$ to denote a universal constant which may take different values from one appearance to another. For simplicity, we take $q_n = (e\sqrt{2\pi}/n^2) \wedge 1$ in (2.14) so that (4.4) holds with $\rho_n = n^{-3}$.

Let $\varepsilon_n$ and $x_*$ be as in Theorem 4 and $\widetilde{L}(\rho) = \sqrt{-\log(2\pi\rho^2)}$ be as in Propositions 1 and 4. With $\rho_n = n^{-3}$, set

$$\eta = \frac{\rho_n}{n} = \frac{1}{n^4}, \qquad \eta^* = \frac{\eta}{\rho_n}\{3\widetilde{L}(\eta)+2\}, \qquad M = \frac{2n\varepsilon_n^2}{(\log n)^{3/2}}. \tag{A.15}$$

Let $x^* = \max(x_*, 1)$ and $\{t^*_{H_j}(\cdot;\rho_n), j \le N\}$ be a $(2\eta^*)$-net of

$$\mathscr{T}_{\rho_n} \cap \{t^*_G : d(f_G, f_{G_n}) \le x^*\varepsilon_n\} \tag{A.16}$$

under the $\|\cdot\|_{\infty,M}$ seminorm as in Proposition 3, with distributions $H_j$ satisfying $d(f_{H_j}, f_{G_n}) \le x^*\varepsilon_n$ and $N = N(\eta^*, \mathscr{T}_{\rho_n}, \|\cdot\|_{\infty,M})$. It is a $(2\eta^*)$-net due to the additional requirements on $H_j$. Since $M \ge 4\sqrt{\log n}$ and $\eta = 1/n^4$ by (4.2) and (A.15), Proposition 3 and (A.15) give

$$\log N \le M_0(\log n)^{3/2} M/2 \le M_0 n\varepsilon_n^2. \tag{A.17}$$

We divide the $\ell_2$ distance of the error into five parts:

$$\|t^*_{\widehat{G}_n}(\mathbf{X};\rho_n) - \boldsymbol{\theta}\| \le \sqrt{nR^*(G_n)} + \sum_{j=1}^4 \zeta_{jn},$$

where $\zeta_{jn}$ are as in (4.9), (4.10), (4.11) and (4.12). As we have mentioned in the outline, the problem is to bound $E_{n,\boldsymbol{\theta}}\zeta_{jn}^2$ in view of (4.13).

Let $A_n$ and $\zeta_{1n}$ be as in (4.9). Since $x^* = 1 \vee x_* \ge 1$ and $n\varepsilon_n^2 \ge 2(\log n)^2$ by (4.2), Theorem 4 gives $P_{n,\boldsymbol{\theta}}\{A_n^c\} \le \exp(-(x^*)^2 n\varepsilon_n^2/(2\log n)) \le 1/n$. Thus, since $\widetilde{L}^2(\rho_n) = -\log(2\pi/n^6)$ with $\rho_n = n^{-3}$, Proposition 1 gives

$$\begin{aligned} E_{n,\boldsymbol{\theta}}\zeta_{1n}^2 &= E_{n,\boldsymbol{\theta}}\sum_{i=1}^n \{(t^*_{\widehat{G}}(X_i;\rho_n) - X_i) + (X_i - \theta_i)\}^2 I_{A_n^c} \\ &\le 2n\widetilde{L}^2(\rho_n)P_{n,\boldsymbol{\theta}}\{A_n^c\} + 2E_{n,\boldsymbol{\theta}}\sum_{i=1}^n (X_i-\theta_i)^2 I_{A_n^c} \\ &\le M_0\log n + 2n\int_0^\infty \min(P\{|N(0,1)|>x\}, 1/n)\,dx^2. \end{aligned}$$

Since $P\{N(0,1)>x\} \le e^{-x^2/2}$ and $\int_0^\infty \min(ne^{-x^2/2},1)\,dx^2/2 = 1+\log n$,

$$E_{n,\boldsymbol{\theta}}\zeta_{1n}^2 \le M_0\log n \le M_0 n\varepsilon_n^2. \tag{A.18}$$

Consider $\zeta_{2n}^2$. Since $t^*_{H_j}(\cdot;\rho_n)$ form a $(2\eta^*)$-net of (A.16) under $\|\cdot\|_{\infty,M}$ and $|t^*_G(x;\rho) - x| \le \widetilde{L}(\rho)$ by Proposition 1, it follows from (4.10) that

$$\begin{aligned} \zeta_{2n}^2 &\le \min_{j\le N}\|t^*_{\widehat{G}_n}(\mathbf{X};\rho_n) - t^*_{H_j}(\mathbf{X};\rho_n)\|^2 I_{A_n} \\ &\le (2\eta^*)^2\#\{i : |X_i|\le M\} + \{2\widetilde{L}(\rho_n)\}^2\#\{i : |X_i| > M\}. \end{aligned}$$

By (4.2), $(n\varepsilon_n^2/\log n)^{p+1} \ge n\{\sqrt{\log n}\,\mu_p^w(G_n)\}^p$, so that by (4.1) and (A.15)

$$\begin{aligned} \int_{|u|\ge M/2} G_n(du) &\le \left(\frac{\mu_p^w(G_n)}{M/2}\right)^p \\ &\le \left(\frac{2n\varepsilon_n^2}{M(\log n)^{3/2}}\right)^p \frac{\varepsilon_n^2}{\log n} = \frac{\varepsilon_n^2}{\log n}. \end{aligned} \tag{A.19}$$

Thus, since $\eta^* = n^{-1}\{3\widetilde{L}(n^{-4}) + 2\}$ and $M \ge 4\sqrt{\log n}$ by (A.15) and (4.2),

$$\begin{aligned} E_{n,\boldsymbol{\theta}}\zeta_{2n}^2 &\le n(2\eta^*)^2 + 4\widetilde{L}^2(n^{-4})E_{n,\boldsymbol{\theta}}\#\{i : |X_i| > M\} \\ &\le M_0(\log n)n\left[\frac{1}{n^2} + \int_{|u|\ge M/2} G_n(du) + P\{|N(0,1)| > 2\sqrt{\log n}\}\right] \\ &\le M_0(\log n)\left(\frac{1}{n} + \frac{n\varepsilon_n^2}{\log n} + \frac{2}{n}\right). \end{aligned}$$

Since $n\varepsilon_n^2 \ge 2(\log n)^2$ by (4.2), we find

$$E_{n,\boldsymbol{\theta}}\zeta_{2n}^2 \le M_0 n\varepsilon_n^2. \tag{A.20}$$

Now, consider $\zeta_{3n}^2$. Since $\widetilde{L}^2(\rho_n) \le M_0 \log n$, it follows from (4.11), Proposition 4 and (A.17) that

$$\begin{aligned} E_{n,\boldsymbol{\theta}}\zeta_{3n}^2 &= \int_0^\infty P_{n,\boldsymbol{\theta}}\{\zeta_{3n} > x\}\,dx^2 \\ &\le \int_0^\infty \min\{1, N\exp(-x^2/(2\widetilde{L}^4(\rho_n)))\}\,dx^2 \qquad \text{(A.21)} \\ &= 2\widetilde{L}^4(\rho_n)(1+\log N) \le M_0(\log n)^2 n\varepsilon_n^2. \end{aligned}$$

For $\zeta_{4n}^2$, it suffices to apply Theorem 3(ii) with $G_0 = G_n$, $G = H_j$, $\rho = \rho_n = n^{-3}$, $x_0 = M/2$ and $\varepsilon_0 = x^*\varepsilon_n \ge d(f_{H_j}, f_{G_n})$, since

$$\zeta_{4n}^2 \le n \max_{j\le N}\{E_{G_n}\{t_{H_j}^*(Y;\rho_n) - \xi\}^2 - R^*(G_n)\} \qquad \text{(A.22)}$$

by (4.12) and (2.5). It follows from (A.19) that the $M_1$ in Theorem 3(ii) is no greater than

$$\frac{\int_{|u|\ge M/2} G_n(du)}{|\log\rho_n|^3(x^*\varepsilon_n)^2} \le \frac{\varepsilon_n^2/\log n}{(\log n)^3\varepsilon_n^2} \le M_0.$$

Since $M = 2n\varepsilon_n^2/(\log n)^{3/2}$ by (A.15) and $n\varepsilon_n^2 \ge 2(\log n)^2$ by (4.2), the $M_2$ in Theorem 3(ii) is no greater than

$$\frac{2(M/2+1)\rho_n}{(\log\rho_n)^2(x^*\varepsilon_n)^2} \le \frac{2(n\varepsilon_n^2/(\log n)^{3/2}+1)/n^3}{(3\log n)^2\varepsilon_n^2} \le \frac{\sqrt{\log n}+1}{n^2(\log n)^4} \le M_0$$

with $\rho_n = n^{-3}$. Thus, by Theorem 3(ii) and (A.22)

$$\zeta_{4n}^2 \le M_0 n|(\log\rho_n)/3|^3\varepsilon_n^2 = M_0 n\varepsilon_n^2(\log n)^3. \qquad \text{(A.23)}$$

Adding (A.18), (A.20), (A.21) and (A.23) together, we have

$$E_{n,\boldsymbol{\theta}}\left(\sum_{j=1}^4 |\zeta_{jn}|\right)^2 \le M_0 n\varepsilon_n^2(\log n)^3.$$

Since $L_n(\widehat{\boldsymbol{\theta}}, \boldsymbol{\theta}) = \|\widehat{\boldsymbol{\theta}} - \boldsymbol{\theta}\|^2/n$, this and (4.13) complete the proof. $\square$

PROOF OF THEOREM 2. As we have mentioned, by (2.34), (2.35) and (2.36), the adaptive minimaxity (2.28) with $\Theta_n = \Theta_{p,C_n,n}$ follows from (2.37). By (4.1) and (2.29), $\mu_p^w(G_{n,\boldsymbol{\theta}}) \le C$ for $\boldsymbol{\theta} \in \Theta_{p,C,n}$, so that by (4.2) and Theorem 5, $\sup_{\boldsymbol{\theta}\in\Theta_{p,C,n}} \widetilde{r}_{n,\boldsymbol{\theta}}(\widehat{t}_n) \le \varepsilon_{p,C,n}(\log n)^{3/2}$ with

$$\varepsilon_{p,C,n}^2 = \max[2\log n, \{nC^p(\log n)^{p/2}\}^{1/(1+p)}](\log n)/n. \qquad \text{(A.24)}$$

Thus, it suffices to verify that for sequences $C_n$ satisfying (2.30),

$$\varepsilon_{p,C_n,n}^2(\log n)^3/J_{p,C_n} \to 0, \qquad \text{(A.25)}$$

where $J_{p,C} = \min\{1, C^{p\wedge 2}\{1 \vee (2\log(1/C^p))\}^{(1-p/2)_+}\}$ as in (2.37).

We consider three cases. For $C_n^{2\wedge p} > e^{-1/2}$, $J_{p,C_n} \ge e^{-1/2}$ and

$$\varepsilon^2_{p,C_n,n}(\log n)^3 = \max\left[\frac{2(\log n)^5}{n}, \{C_n(\log n)^{9/2+4/p}/n\}^{p/(1+p)}\right] = o(1),$$

since $\kappa_2(p) = 9/2 + 4/p$ in (2.30).

For $p < 2$ and $C_n^p \le e^{-1/2}$, $J_{p,C_n} = C_n^p\{2\log(1/C_n^p)\}^{1-p/2}$, so that by (A.24)

$$\begin{aligned}\varepsilon^2_{p,C_n,n}&(\log n)^3/J_{p,C_n}\\ &= \max\left[\frac{2(\log n)^5}{nC_n^p\{\log(1/C_n^p)\}^{1-p/2}}, \frac{(\log n)^{4+p/(2+2p)}}{(nC_n^p)^{p/(1+p)}\{\log(1/C_n^p)\}^{1-p/2}}\right].\end{aligned}$$

Since the case $C_n^p > n^{-1/2}$ is trivial, it suffices to consider the case $C_n^p \le n^{-1/2}$ where

$$\frac{\varepsilon^2_{p,C_n,n}(\log n)^3}{J_{p,C_n}} \asymp \max\left[\frac{(\log n)^{4+p/2}}{nC_n^p}, \frac{(\log n)^{3+p/2+p/(2+2p)}}{(nC_n^p)^{p/(1+p)}}\right].$$

Since $4 + p/2 \le p\kappa_1(p) = 4 + 3/p + p/2 = (1 + 1/p)\{3 + p/2 + p/(2 + 2p)\}$, (2.30) still implies (A.25).

Finally, for $p \ge 2$ and $C_n^2 \le e^{-1/2}$, $J_{p,C_n} = C_n^2$, so that

$$\frac{\varepsilon^2_{p,C_n,n}(\log n)^3}{J_{p,C_n}} = \max\left[\frac{2(\log n)^5}{nC_n^2}, \left\{\frac{C_n(\log n)^{9/2+4/p}}{nC_n^{2(1+1/p)}}\right\}^{p/(1+p)}\right].$$

Since $nC_n^{1+2/p} = n^{1/2-1/p}(nC_n^2)^{1/2+1/p}$, we need $(\log n)^5/(nC_n^2) \to 0$ for $p > 2$ and $(\log n)^{13/2}/(nC_n^2) \to 0$ for $p = 2$. Again (2.30) implies (A.25). $\square$

Proof of Theorem 7. Since the oracle inequality (4.8) is based on the weak $\ell_p$ norm, the proof of Theorem 2 also provides (4.18). $\square$

Proof of Proposition 5. Let $\widehat{G}_n^*$ be the exact generalized MLE as in (2.12). Since $\varphi(x)$ is decreasing in $|x|$, we have $\widehat{G}_n^*([u_1, u_m]) = 1$. Let $I_j = (u_{j-1}, u_j]$ and $I_j^* = [u_{j-1}, u_j]$ for $j \ge 2$ and $I_1 = I_1^* = \{u_1\}$. Let $H_{m,j}$ be sub-distributions with support $\{u_{j-1}, u_j\} \cap I_j^*$ such that

$$H_{m,j}(I_j^*) = \widehat{G}_n^*(I_j), \qquad \int_{I_j^*} uH_{m,j}(du) = \int_{I_j} u\widehat{G}_n^*(du), \tag{A.26}$$
$$1 \le j \le m.$$

Let $j > 1$ and $x \in [u_1, u_m]$ be fixed. Set $x_j = x - (u_j + u_{j-1})/2$ and $t = u - (u_j + u_{j-1})/2$ for $u \in I_j^*$. Since $|x_j t| \le (u_m - u_1)\varepsilon/2 \le n^{-1/2} \le 1$,

$$-(1 - e^{-t^2/2})e^{x_j t} \le e^{x_j t - t^2/2} - (1 + x_j t) \le x_j^2 t^2 e^{x_j t - t^2/2}, \tag{A.27}$$

where the second inequality follows from $e^{-t^2/2}(1-x_j t) \le e^{-x_j t}$. Since $\varphi(x-u) = \varphi(x_j - t) = \varphi(x_j)\exp(x_j t - t^2/2)$, (A.26) and (A.27) yield

$$\begin{aligned}
&\int_{I_j} \varphi(x-u)\widehat{G}_n^*(du) - \int_{I_j^*} \varphi(x-u)H_{m,j}(du) \\
&\qquad \le \int_{I_j} x_j^2 t^2 \varphi(x-u)\widehat{G}_n^*(du) + \int_{I_j^*} (e^{t^2/2}-1)\varphi(x-u)H_{m,j}(du) \\
&\qquad \le (u_m - u_1)^2(\varepsilon/2)^2 \int_{I_j} \varphi(x-u)\widehat{G}_n^*(du) \\
&\qquad\quad + (e^{\varepsilon^2/8}-1)\int_{I_j^*} \varphi(x-u)H_{m,j}(du).
\end{aligned}$$

Let $H_m = \sum_{j=1}^m H_{m,j}$. Summing the above inequality over $j$, we find $e^{\varepsilon^2/8} \times f_{H_m}(x) \ge (1-\eta) f_{\widehat{G}_n^*}(x)$ with $\eta = \varepsilon^2(u_m - u_1)^2/4 \le 1/n - \varepsilon^2/8$. Thus,

$$\prod_{i=1}^n \frac{f_{H_m}(X_i)}{f_{\widehat{G}_n^*}(X_i)} \ge e^{-n\varepsilon^2/8}(1-\eta)^n \ge e^{-n(\varepsilon^2/8+\eta)} \ge e^{-1}. \tag{A.28}$$

Let $\mathscr{H}_m$ be the set of all distributions with support $\{u_1, \dots, u_m\}$ and $\widehat{G}_n = \sum_{j=1}^m \widehat{w}_j^{(k)} \delta_{u_j}$. The upper bound in [10, 35] and (6.3) provide

$$\sup_{H \in \mathscr{H}_m} \prod_{i=1}^n \frac{f_H(X_i)}{f_{\widehat{G}_n}(X_i)} \le \max_{j \le m} \left( \frac{w_j^{(k)}}{w_j^{(k-1)}} \right)^n \le \frac{1}{e q_n}.$$

This and (A.28) imply $\prod_{i=1}^n f_{\widehat{G}_n^*}(X_i) \le q_n^{-1} \prod_{i=1}^n f_{\widehat{G}_n}(X_i)$. $\square$

Department of Statistics
and Biostatistics
Hill Center
Busch Campus
Rutgers University
Piscataway, New Jersey 08854
E-mail: wenhua@eden.rutgers.edu
czhang@stat.rutgers.edu